\documentclass[a4paper,12pt]{article}
\usepackage{amsmath,amsthm,amsfonts,amssymb,bm,mathrsfs}
\usepackage{enumerate}
\usepackage[numbers]{natbib}
\usepackage{natbib}
\usepackage{xcolor}
\usepackage{sgame}

\usepackage{setspace}
\renewcommand{\baselinestretch}{0.95}

\newtheorem{defn}{Definition}
\newtheorem{thm}{Theorem}
\newtheorem{lem}{Lemma}
\newtheorem*{claim}{Claim}

\newtheorem{prop}{Proposition}
\newtheorem{rmk}{Remark}

\newtheorem{coro}{Corollary}

\newcommand{\bN}{\mathbb{N}}
\newcommand{\bR}{\mathbb{R}}
\newcommand{\cB}{\mathcal{B}}

\newcommand{\cF}{\mathcal{F}}
\newcommand{\cG}{\mathcal{G}}

\newcommand{\cI}{\mathcal{I}}

\newcommand{\cM}{\mathcal{M}}

\newcommand{\cR}{\mathscr{R}}

\newcommand{\cT}{\mathcal{T}}
\newcommand{\cU}{\mathcal{U}}

\DeclareMathOperator{\rmd}{d\!}

\usepackage{hyperref}

\hypersetup{colorlinks=true,
            citecolor=blue,
            linkcolor=blue,
            urlcolor=blue,
            bookmarksopen=true,
            pdfstartview={XYZ null null 1.00},
            pdfpagelayout={SinglePage}
}

\allowdisplaybreaks

\begin{document}
\title{Conditional Expectation of Banach Valued Correspondences\footnote{The main results in Sections 3, 4 and 5 of  this paper were included in Chapter~3 of Wei He's 2014 Ph.D dissertation (\citet{He2014}).}}
\author{Wei~He\thanks{Department of Economics, The Chinese University of Hong Kong, Shatin, N.T., Hong Kong SAR. E-mail: hewei@cuhk.edu.hk.}
\and
Yeneng~Sun\thanks{Departments of Economics and Mathematics, National University of Singapore, 10 Lower Kent Ridge Road, Singapore 119076. E-mail: ynsun@nus.edu.sg.}
}
\date{\today}
\maketitle

\abstract{We present some regularity properties (convexity, weak/weak$^*$ compactness and preservation of weak/weak$^*$ upper hemicontinuity) for Bochner/Gel$^\prime$fand conditional expectation of Banach valued correspondences under the nowhere equivalence condition. These regularity properties for Bochner/Gel$^\prime$fand integral of Banach valued correspondences are obtained as corollaries. Similar properties for regular conditional distributions are also covered by the corresponding results for Gel$^\prime$fand conditional expectation of correspondences. We prove the necessity of the nowhere equivalence condition for any of these properties to hold. As an application, we show that the nowhere equivalence condition is satisfied on the space of players if and only if a pure-strategy Nash equilibrium exists in a general class of large games.

%\textbf{JEL classification}:

}

%\newpage
%
%\tableofcontents

\newpage

\section{Introduction}\label{sec-intro}

The theory of correspondences has found wide applications in a number of areas including optimization, control theory, stochastic analysis, game theory and general equilibrium theory.\footnote{See, for example, the books \citet{AB2006}, \citet{AF1990}, \citet{Hd1974}, \citet{Mi}, and also the papers \citet{KS2014}, \citet{KZ2012}, \citet{Podczeck2008}, \citet{SY2008}, \citet{SZ2015} and \citet{Yu2013} for some recent developments and applications.} The most useful regularity properties for the integral of correspondences over an atomless measure space are convexity, compactness and preservation of upper hemicontinuity.
%, since they are the key properties when one needs to apply the fixed-point argument.
%\footnote{The integral theory of correspondences has been well-developed. See, for example, \citet{AB2006} and \citet{Hd1974} for results in the finite-dimensional setting, and \citet{KS1996}, \citet{Sun1997}, \citet{Podczeck2008} and \citet{SY2008} for results in the infinite-dimensional setting.}
Given the relationship between conditional expectation and integration, it is natural to consider the conditional expectation of correspondences, and expect that such a generalization could find more applications in various fields. The main aim of this paper is to characterize the regularity properties of convexity, compactness and preservation of upper hemicontinuity for Bochner/Gel$^\prime$fand conditional expectation of Banach valued correspondences.\footnote{As noted in Remark~\ref{rmk-distribution} and Corollary~\ref{coro-convex distribution} below, such properties for the case of Gel$^\prime$fand conditional expectation imply the same properties for regular conditional distribution (which is a standard concept in probability).}

%The main aim of this paper is to present a comprehensive theory of conditional expectation of correspondences in the infinite-dimensional setting.

When the range space is finite dimensional, the theory on the integral (or more generally, conditional expectation) of correspondences has been considered in the literature. However, it is well known that even if one considers the integral of correspondences, all the regularity properties (convexity, compactness, and preservation of upper hemicontinuity) may fail when the range space is an infinite-dimensional Banach space.\footnote{It is well known that those regularity properties hold for integral of correspondences in the finite-dimensional setting; see, for example, \citet{AB2006} and \citet{Hd1974}. Such results are based on Lyapunov's theorem, which has been extensively studied in the literature; for some recent developments of Lyapunov's theorem, see, for example, \citet{DF2011, DF2014} and \citet{FP2010}. \citet{Sole1970} presented a conditional version of this theorem in the finite-dimensional setting. The conditional expectation of correspondences in the finite-dimensional setting has been studied systematically in \citet{HS2018b}. For the infinite-dimensional setting, see \citet{DU1977} and \citet{Sun1997} for some counterexamples on the failure of those regularity properties, while \citet{KS1996} showed that those properties hold for integral of correspondences with a countable range. \label{fn-main}} Thus, these regularity properties would also fail in general for conditional expectation of Banach valued correspondences.

Rather than working with the conditional expectation of Banach valued correspondences as in this paper, the earlier works consider its closure; see, for example, \citet{HU1977}, \citet{Hiai1985} and \citet{Papageorgiou1987}. As noted in the paragraph above and Footnote~\ref{fn-main}, the regularity properties as considered here fail in general for the conditional expectation of correspondences without taking the closure. To resolve this issue, we work with the condition of ``nowhere equivalence'', which captures in a certain sense relative richness between a $\sigma$-algebra and its sub-$\sigma$-algebra. Based on the nowhere equivalence condition, we are able to characterize the desirable properties for both Bochner and Gel$^\prime$fand conditional expectation of Banach valued correspondences.

Formally, let $(T, \cT, \lambda)$ be an atomless probability space and $\cF$ a countably generated sub-$\sigma$-algebra of $\cT$. Then $\cT$ is said to be nowhere equivalent to $\cF$ if these two $\sigma$-algebras do not coincide (modulo null sets) on any non-negligible subset of $T$. In Theorems~\ref{thm-ce-bochner} and \ref{thm-ce-G}, we show that if a Banach valued correspondence $G$ is $\cF$-measurable and its selections are allowed to be $\cT$-measurable, then all of these regularity properties hold for the Bochner and Gel$^\prime$fand conditional expectation of the correspondence $G$, respectively.

Theorems~\ref{thm-ce-bochner} and \ref{thm-ce-G} indicate that the nowhere equivalence condition is sufficient for the relevant regularity properties to hold.  We show in Theorem~\ref{thm-converse} that the nowhere equivalence condition is also necessary for the validity of any of those regularity properties. Thus, the nowhere equivalence condition is ``minimal'' if one needs to work with any of those properties for conditional expectation of Banach valued correspondences in general.

It is clear that the integral can be viewed as the conditional expectation with respect to the trivial $\sigma$-algebra, and any atomless $\sigma$-algebra is obviously nowhere equivalent to the trivial $\sigma$-algebra. Thus, the results on the Bochner/Gel$^\prime$fand integral of Banach valued correspondences can be delivered as direct corollaries (Proposition~\ref{prop-Bochener}/Proposition~\ref{prop-G integral}). Furthermore, the norm compactness and preservation of norm upper hemicontinuity for integral of Banach valued correspondences are also presented in Proposition~\ref{prop-norm integral}. Our Propositions~\ref{prop-Bochener}-\ref{prop-G integral} generalize the results in a few earlier papers on integral of Banach valued correspondences by using the more general condition of nowhere equivalence.\footnote{To obtain convexity for Bochner integral of Banach valued correspondences, \citet{RY1991} proposed to work with measure spaces with the associated $L^{\infty}$ spaces over any non-null measurable set having strictly larger cardinality than that of the continuum. 
\citet{Sun1997} resolved the issue on the regularity properties 
(convexity, weak/weak$^*$ compactness
and preservation of weak/weak$^*$ upper hemicontinuity)
for Bochner and Gel$^\prime$fand  integral of correspondences based on atomless Loeb spaces; for the construction of a Loeb space, see \citet{Loeb1975}. These results were further generalized by \citet{SY2008} and \citet{Podczeck2008} to those over a saturated probability space. Note that a probability space is saturated if its $\sigma$-algebra is nowhere equivalent to any countably-generated sub-$\sigma$-algebra; see \citet{KS2009}.}  %These papers did not consider the conditional expectation of correspondences.
%Propositions~\ref{prop-Bochener}-\ref{prop-G integral} extend the results in those works by using the more general condition of nowhere equivalence.

%the results in \citet{Sun1997}, \citet{SY2008} and \citet{Podczeck2008}.

%Our Propositions~\ref{prop-Bochener}-\ref{prop-G integral} on the regularity properties for Bochner and Gel$^\prime$fand  integral of correspondences follow from Theorems~\ref{thm-ce-bochner} and \ref{thm-ce-G} easily, which generalize the results in \citet{Sun1997}, \citet{SY2008} and \citet{Podczeck2008}.\footnote{\citet{Sun1997} resolved the issue on regularity properties for Bochner and Gel$^\prime$fand  integral of correspondences based on atomless Loeb spaces; for the construction of a Loeb space, see \citet{Loeb1975}. These results were further generalized by \citet{SY2008} and \citet{Podczeck2008} to those over a saturated probability space. Note that a probability space is saturated if its $\sigma$-algebra is nowhere equivalent to any countably-generated sub-$\sigma$-algebra; see \citet{KS2009}.}  These papers did not consider the conditional expectation of correspondences. Propositions~\ref{prop-Bochener}-\ref{prop-G integral} extend the results in these papers, using the more general condition of nowhere equivalence.

To illustrate the usefulness of the characterization results, we apply Theorems~\ref{thm-ce-bochner} and \ref{thm-ce-G} to study a new class of large games and obtain new equilibrium existence results. In a typical large game, each player's payoff depends on her action and the societal aggregate, which is often modeled as the average actions from all the players.\footnote{Large games of this type are also called mean-field games; see, for example, \citet{CHM2017}, \citet{CL2015}, \citet{HCM2012}, \citet{Lacker2020}, \citet{LL2007}, \citet{Nutz2018} and \citet{Nutz-Zhang2019}.} Such modeling of ``aggregate actions'' completely removes the differences among individual players. As argued in \citet{KRSY2013}, economic agents are often associated with biological or social traits, and hence it is possible that an individual player's payoff depends on the aggregate actions in a rich manner. In Section~\ref{sec-game}, we consider a large game where the action space is a subset of some infinite-dimensional Banach space, and the externality part is the conditional expectation of players' action profiles, given their characteristics. That is, a player may care about the aggregate actions from those players who share the same characteristics.\footnote{For example, one may care about not only the aggregate actions of the whole society, but also the aggregate actions from each group with the same income level.} We show in Theorem~\ref{thm-existence} that the nowhere equivalence condition characterizes the pure-strategy equilibrium existence in such large games in the sense that a pure-strategy Nash equilibrium exists if and only if the nowhere equivalence condition is satisfied on the player space.\footnote{In the end of Section~\ref{sec-game}, we provide more discussions on the equilibrium existence results in the literature, and compare those results with our Theorem~\ref{thm-existence}.}

The paper is organized as follows. Section~\ref{sec-basic} collects some basic definitions of correspondences and also the nowhere equivalence condition. In Sections~\ref{sec-Bochner} and \ref{sec-G integral}, we prove the regularity properties for Bochner and Gel$^\prime$fand conditional expectation of correspondences, respectively. Section~\ref{sec-converse} presents the necessity result. As an application of main results, Section~\ref{sec-game} characterizes the existence of pure-strategy Nash equilibrium via the nowhere equivalence condition. All the proofs are collected in Appendix.

%%%%%%%%%%%%%%%%%%%%%%%%%%%%%%%%%%%%%%%%%%%%%%%%%%%%%%%%%%%%%%%%%%%%%%%%%%%%%%%%%%%%%%%%%%%%%%%%%%%%%%%%%%%%%%%%%%%%%%%%%%%%%%%%%%%%%%%%%%%%%%%%%%%%%%%%%%%%%%%%%%%%%%%%%%%%%%%%%%%%%%%%%%

\section{Basics}\label{sec-basic}

Let $X$ be a Banach space endowed with the norm $\|\cdot\|$, and $(T,\cT,\lambda)$ a complete probability space endowed with a countably additive probability measure $\lambda$.

Given a mapping $f$ from $(T,\cT,\lambda)$ to $X$, if $f$ can be approximated in norm by a sequence of simple functions and $\|f\|$ is integrable on $(T,\cT,\lambda)$, then $f$ is said to be Bochner integrable.

A mapping $f$ from $(T, \cT, \lambda)$ to the norm dual $X^*$ of a Banach space $X$ is said to be Gel$^\prime$fand integrable if $f(\cdot)(x)$ is integrable over $(T, \cT, \lambda)$ for each $x\in X$. The Gel$^\prime$fand integral of $f$ is the unique element $x^*$ in $X^*$ such that $x^*(x) = \int_T f(\cdot)(x) \rmd\lambda$ for all $x\in X$, and $x^*$ is denoted by $\int_T f \rmd\lambda$.

For a Bochner (resp. Gel$^\prime$fand) integrable mapping $f$ from $T$ to some Banach space $Y$, a mapping $g$ is called a Bochner (resp. Gel$^\prime$fand) conditional expectation of $f$ given a sub-$\sigma$-algebra $\cF \subseteq \cT$, if $g$ is $\cF$-measurable and $\int_D f \rmd \lambda = \int_D g \rmd \lambda$ for any $D \in \cF$. Hereafter, we denote the conditional expectation of $f$ given $\cF \subseteq \cT$ by $E(f|\cF)$.

Given two nonempty sets $T$ and $Y$, a correspondence $F$ from $T$ to $Y$ is a mapping from the set $T$ to the collection of all nonempty subsets of $Y$. Suppose that $Y$ is a topological space endowed with the Borel $\sigma$-algebra $\cB(Y)$. A correspondence $F$ from $(T,\cT,\lambda)$ to $Y$ is said to be measurable if $\{t\in T\colon F(t)\cap E \neq \emptyset\} \in \cT$ for each closed subset $E\subseteq Y$. A correspondence $F$ from $(T,\cT,\lambda)$ to $Y$ is said to be convex (resp. closed, compact) valued if $F(t)$ is convex (resp. closed, compact) for $\lambda$-almost all $t\in T$. A mapping $f$ is called a selection of $F$ if $f(t) \in F(t)$ for $\lambda$-almost all $t \in T$.

For $1\le p \le \infty$, a correspondence $F$ is said to be $p$-integrably bounded if there is a real-valued function $h$ such that $h^p$ is integrable on $(T,\cT,\lambda)$, and  $\sup\{\|x\|: x\in F(t)\} \leq h(t)$ for $\lambda$-almost all $t\in T$. If $p=1$, then we just say that $F$ is integrably bounded.

A correspondence $F$ from a topological space $Y$ to another  topological space $Z$ is said to be upper hemicontinuous at $y_0\in Y$, if for any open set $O_Z$ that contains $F(y_0)$, there exists an open neighborhood $O_Y$ of $y_0$ such that for any $y\in O_Y$, $F(y)\subseteq O_Z$. The correspondence $F$ is said to be upper hemicontinuous if it is upper hemicontinuous at every point $y\in Y$.

%The usual Lebesgue unit interval is denoted by $(I,\cB,\eta)$, \textit{i.e.}, the unit interval $I=[0,1]$ endowed with the Borel $\sigma$-algebra $\cB$ and the Lebesgue measure $\eta$.

Hereafter, we assume that $(T,\cT,\lambda)$ is an atomless probability space. Given a nonnegligible subset $D \in \cT$, we shall define the restricted space $(D,\cT^D,\lambda^D)$ of $(T,\cT,\lambda)$ on $D$ as follows: let $\cT^D$ be the $\sigma$-algebra $\{ D\cap D' \colon D' \in \cT \}$, and $\lambda^D$ the probability measure re-scaled from the restriction of $\lambda$ to $\cT^D$. For a sub-$\sigma$-algebra $\cF \subseteq \cT$, one can define the restricted space $(D,\cF^D,\lambda^D)$ similarly. Hereafter, we assume that $\cF$ is countably generated. The following condition compares the $\sigma$-algebra $\cT$ with its sub-$\sigma$-algebra $\cF$.

\begin{defn}\label{defn-nq}
The $\sigma$-algebra $\cT$ is said to be nowhere equivalent to its sub-$\sigma$-algebra $\cF$ if for every $D\in\cT$ with $\lambda(D)>0$, there exists a $\cT$-measurable set $D_0 \subseteq D$ such that $\lambda(D_0\triangle D_1)>0$ for any $D_1\in\cF^D$.\footnote{For more discussions on the condition of nowhere equivalence, see \citet{HS2018a}.}
\end{defn}

%%%%%%%%%%%%%%%%%%%%%%%%%%%%%%%%%%%%%%%%%%%%%%%%%%%%%%%%%%%%%%%%%%%%%%%%%%%%%%%%%%%%%%%%%%%%%%%%%%%%%%%%%%%%%%%%%%%%%%%%%%%%%%%%%%%%%%%%%%%%%%%%%%%%%%%%%%%%%%%%%%%%%%%%%%%%%%%%%%%%%%%%%%

\section{Bochner integral and conditional expectation of correspondences}\label{sec-Bochner}

In this section, we consider the Bochner integral and conditional expectation of Banach valued correspondences.\footnote{As we always work with Banach valued correspondences, hereafter we often just say integral/conditional expectation of correspondences for simplicity.} Specifically, we shall show that the condition of nowhere equivalence is sufficient for several desirable regularity properties (convexity, compactness, and preservation of upper hemicontinuity) of Bochner integral and conditional expectation of correspondences. In Section~\ref{sec-converse}, we shall show that the nowhere equivalence condition is also necessary for these properties.

Fix the Banach space $X$ and the atomless probability space $(T,\cT,\lambda)$. Let $\cG$ be a sub-$\sigma$-algebra of $\cT$. Given a correspondence $F$ from $T$ to $X$, denote
$$I_{F}^{\cG} = \bigg\{\int_T f(t) \rmd\lambda(t):f\text{ is a } \cG\text{-measurable Bochner integrable selection of }F \bigg\},
$$
and
$${CI}_{F}^{(\cT,\cG)} = \bigg\{ E(f|\cG):f\mbox{ is a } \cT\mbox{-measurable Bochner integrable selection of }F \bigg\}.
$$
For $p\geq 1$, let $L^{\cG}_p(T,X)$ be the set of all $\cG$-measurable Bochner integrable mappings from $T$ to $X$ under the usual norm; that is, for $1\leq p <\infty$,
$$L^{\cG}_p(T,X) = \bigg\{f\colon f \mbox{ is } \cG \mbox{-measurable and } \left(\int_T \|f\|^p \rmd\lambda \right)^{\frac{1}{p} } < \infty \bigg\}
$$
and
$$L^{\cG}_{\infty}(T, X) = \bigg\{f\colon f \mbox{ is } \cG \mbox{-measurable, } \exists \, c > 0 \mbox{ such that }\lambda\left(t\colon \|f(t)\| < c  \right)=1 \bigg\}.
$$

In this section, we assume that $X$ is a separable Banach space, and $X^*$ has the Radon-Nikod\'ym property.  Since $X$ is separable, it is straightforward to check that $L^{\cG}_p(T,X)$ is also a separable Banach space for $1 \leq p <\infty$ when $\cG$ is countably generated. As $X^*$ has the Radon-Nikod\'ym property,  by Theorem 1 in \citet[p.98]{DU1977}, $L^{\cG}_{q}(T, X^*)$ consists of the continuous linear functionals of $L_{p}^\cG(T,X)$  for $1 \leq p < \infty$, where $X^*$ is the dual space of $X$ and $\frac{1}{p}+\frac{1}{q}=1$.

We shall start by considering the properties on the convexity, compactness, and preservation of upper hemicontinuity for conditional expectation of Banach-valued correspondences. The corresponding properties for the Bochner integral of correspondences follow directly.

\begin{thm}\label{thm-ce-bochner}
If $\cT$ is nowhere equivalent to $\cF$, then for any sub-$\sigma$-algebra $\cG$ of $\cF$, we have the following properties.

\begin{enumerate}
  \item[A1] For any $\cF$-measurable and weakly compact valued correspondence $F$ from $T$ to $X$, ${CI}_{F}^{(\cT,\cG)}$ is convex.

  \item[A2] For any $\cF$-measurable, $p$-integrably bounded and weakly compact valued correspondence $F$ from $T$ to $X$ with $1\leq p < \infty$, $CI_{F}^{(\cT,\cG)}$ is weakly compact in $L_p^\cG(T,X)$.

  \item[A3] For any $\cF$-measurable, $p$-integrably bounded and weakly compact valued correspondence $F$ from $T$ to $X$  with $1\leq p < \infty$,
  $$CI_{F}^{(\cT,\cG)}=CI_{\overline{co}F}^{(\cT,\cG)},
  $$
  where $\overline{co}F$ is the correspondence such that $\overline{co}F(t)$ is the norm closure of the convex hull of $F(t)$ for each $t\in T$.

  \item[A4] Suppose that $G$ is an $\cF$-measurable, $p$-integrably bounded and weakly compact valued correspondence from $T$ to $X$ with $1 \le p < \infty$. Let $F \colon T\times Y\to X$ be a weakly closed valued correspondence ($Y$ is a metric space). If
      \begin{enumerate}
        \item $F(t,y)\subseteq G(t)$ for $\lambda$-almost all $t\in T$;
        \item for any $y\in Y$, $F(\text{\ensuremath{\cdot}},y)$ is $\cF$-measurable;
        \item for any $t\in T$, $F(t,\cdot)$ is weakly upper hemicontinuous;
      \end{enumerate}
      then $H(y)= CI_{F_y}^{(\cT,\cG)}$ is weakly upper hemicontinuous in $L_p^\cG(T,X)$.
\end{enumerate}
\end{thm}

\citet{HU1977}, \citet{Hiai1985} and \citet{Papageorgiou1987} studied the properties for the closure of ${CI}_{F}^{(\cT,\cG)}$; that is, the closure of the Bochner conditional expectation of the correspondence $F$. In contrast, we directly show the regularity properties for ${CI}_{F}^{(\cT,\cG)}$ rather than its closure.

The regularity properties for integral of correspondences follow from those for conditional expectation of correspondences. In particular, if we take the sub-$\sigma$-algebra $\cG$ to be the trivial $\sigma$-algebra $\{T, \emptyset\}$ in Theorem~\ref{thm-ce-bochner}, then the conditional expectation is reduced to be the infinite-dimensional integral, and the following result is immediate.

\begin{prop}\label{prop-Bochener}
If $\cT$ is nowhere equivalent to $\cF$, then we have the following properties.

\begin{enumerate}
  \item[B1] For any $\cF$-measurable and weakly compact valued correspondence $F$ from $T$ to $X$, $I_{F}^{\cT}$ is convex.

  \item[B2] For any $\cF$-measurable, integrably bounded and weakly compact valued correspondence $F$ from $T$ to $X$, $I_{F}^{\cT}$ is weakly compact.

  \item[B3] For any $\cF$-measurable, integrably bounded and weakly compact valued correspondence $F$ from $T$ to $X$, $I_{F}^{\cT}=I_{\overline{co}F}^{\cT}$.

  \item[B4] Suppose that $G$ is an $\cF$-measurable, integrably bounded and weakly compact valued correspondence from $T$ to $X$. Let $F \colon T\times Y\to X$ be a weakly closed valued correspondence ($Y$ is a metric space). If
      \begin{enumerate}
        \item $F(t,y)\subseteq G(t)$ for $\lambda$-almost all $t\in T$;
        \item for any $y\in Y$, $F(\text{\ensuremath{\cdot}},y)$ is $\cF$-measurable;
        \item for any $t\in T$, $F(t,\cdot)$ is weakly upper hemicontinuous;
      \end{enumerate}
      then $H(y)=I_{F_{y}}^{\cT}$ is weakly upper hemicontinuous.
\end{enumerate}
\end{prop}

In the following proposition, we show that the properties on norm compactness and preservation of norm upper hemicontinuity of integral of correspondences also hold.

\begin{prop}\label{prop-norm integral}
If $\cT$ is nowhere equivalent to $\cF$, then we have the following properties.

\begin{enumerate}
  \item[B5] For any $\cF$-measurable, integrably bounded and norm compact valued correspondence $F$ from $T$ to $X$, $I_{F}^{\cT}$ is norm compact.

  \item[B6] Suppose that $G$ is an $\cF$-measurable, integrably bounded and norm compact valued correspondence from $T$ to $X$. Let $F \colon T\times Y\to X$ be a weakly closed valued correspondence ($Y$ is a metric space). If
      \begin{enumerate}
        \item $F(t,y)\subseteq G(t)$ for $\lambda$-almost all $t\in T$;
        \item for any $y\in Y$, $F(\text{\ensuremath{\cdot}},y)$ is $\cF$-measurable;
        \item for any $t\in T$, $F(t,\cdot)$ is norm upper hemicontinuous;
      \end{enumerate}
      then $H(y)=I_{F_{y}}^{\cT}$ is norm upper hemicontinuous.
\end{enumerate}
\end{prop}

\begin{rmk}
The condition of the Radon-Nikod\'ym property on $X^*$ is not needed in Propositions~\ref{prop-Bochener} and \ref{prop-norm integral}. In the case that $\cG$ is non-trivial,  this condition is to guarantee that $L^{\cG}_{q}(T, X^*)$ consists of all the continuous linear functionals of $L_{p}^\cG(T,X)$ for $1 \leq p < \infty$ with $\frac{1}{p}+\frac{1}{q}=1$; that is, $L^{\cG}_{q}(T, X^*)$ is the dual space of $L_{p}^\cG(T,X)$. When $\cG$ is the trivial $\sigma$-algebra $\{T, \emptyset\}$, $L^{\cG}_{q}(T, X^*)$ and $L_{p}^\cG(T,X)$ are simply $X^*$ and $X$. As $X^*$ is the dual space of $X$, the condition of Radon-Nikod\'ym property on $X^*$ is unnecessary.
\end{rmk}

%%%%%%%%%%%%%%%%%%%%%%%%%%%%%%%%%%%%%%%%%%%%%%%%%%%%%%%%%%%%%%%%%%%%%%%%%%%%%%%%%%%%%%%%%%%%%%%%%%%%%%%%%%%%%%%%%%%%%%%%%%%%%%%%%%%%%%%%%%%%%%%%%%%%%%%%%%%%%%%%%%%%%%%%%%%%%%%%%%%%%%%%%%

\section{Gel$^\prime$fand integral and conditional expectation of correspondences}\label{sec-G integral}

In this section, we turn to the Gel$^\prime$fand integral and conditional expectation of correspondences. As in Section~\ref{sec-Bochner}, we shall show that the condition of nowhere equivalence is sufficient for the desirable regularity properties (convexity, compactness, and preservation of upper hemicontinuity) of Gel$^\prime$fand integral and conditional expectation of correspondences.

Throughout this section, we assume that $X$ is a separable Banach space. Let $\cG$ be a sub-$\sigma$-algebra of $\cT$. Given a correspondence $F$ from $T$ to  $X^*$,  denote
$$G_{F}^{\cG} = \left\{ \int_T f \rmd\lambda(t):f\text{ is a } \cG\text{-measurable Gel$^\prime$fand integrable selection of }F \right\},
$$
and
$${CG}_{F}^{(\cT,\cG)} = \left\{ E(f|\cG):f\mbox{ is a } \cT\mbox{-measurable Gel$^\prime$fand  integrable selection of }F \right\},
$$
where the integral and the conditional expectation take the Gel$^\prime$fand integral.

\begin{thm}\label{thm-ce-G}
If $\cT$ is nowhere equivalent to $\cF$, then for any sub-$\sigma$-algebra $\cG$ of $\cF$, we have the following properties.

\begin{enumerate}
  \item[C1] For any $\cF$-measurable and weak$^*$ compact valued correspondence $F$ from $T$ to $X^*$, ${CG}_{F}^{(\cT,\cG)}$ is convex.

  \item[C2] For any $\cF$-measurable, $p$-integrably bounded, and weak$^*$~compact valued correspondence $F$ from $T$ to $X^*$ with $1 < p \leq \infty$, $CG_{F}^{(\cT,\cG)}$ is weak$^*$~compact (in the dual space of $L_q^\cG(T,X)$ with $\frac{1}{p} + \frac{1}{q} = 1$).

  \item[C3] For any $\cF$-measurable, $p$-integrably bounded, and weak$^*$~compact valued correspondence $F$ from $T$ to $X^*$ with $1 < p \leq \infty$, $CG_{F}^{(\cT,\cG)} = CG_{w^*-\overline{co}F}^{(\cT,\cG)}$, where $w^*-\overline{co}F$ is the correspondence such that for each $t\in T$, $w^*-\overline{co}F(t)$ is the weak$^*$ closure of the convex hull of $F(t)$.

  \item[C4] Suppose that $G$ is an $\cF$-measurable, $p$-integrably bounded, and weak$^*$ compact valued correspondence from $T$ to $X^*$ with $1< p \leq \infty$. Let $F$ be a weak$^*$ closed valued correspondence from $T \times Y$ to $X^*$ ($Y$ is a metric space). If
      \begin{enumerate}
        \item $F(t,y) \subseteq G(t)$ for $\lambda$-almost all $t\in T$;
        \item for any $y\in Y$, $F(\text{\ensuremath{\cdot}},y)$ is $\cF$-measurable;
        \item for any $t\in T$, $F(t,\cdot)$ is weak$^*$ upper hemicontinuous;
      \end{enumerate}
      then $H(y)= CG_{F_y}^{(\cT,\cG)}$ is weak$^*$ upper hemicontinuous (in the dual space of $L_q^\cG(T,X)$).
\end{enumerate}
\end{thm}

Proposition~\ref{prop-G integral} below considers the Gel$^\prime$fand integral of correspondences. In Theorem~\ref{thm-ce-G}, if we take the sub-$\sigma$-algebra $\cG$ to be the trivial $\sigma$-algebra $\{T, \emptyset\}$, then the Gel$^\prime$fand conditional expectation is reduced to be the Gel$^\prime$fand integral.

\begin{prop}\label{prop-G integral}
If $\cT$ is nowhere equivalent to $\cF$, then we have the following properties.

\begin{enumerate}
  \item[D1] For any $\cF$-measurable and weak$^*$ compact valued correspondence $F$ from $T$ to $X^*$, $G_{F}^{\cT}$ is convex.

  \item[D2] For any $\cF$-measurable, $p$-integrably bounded, and weak$^*$~compact valued correspondence $F$ from $T$ to $X^*$ with $1 < p \leq \infty$,  $G_{F}^{\cT}$ is weak$^*$~compact in $X^*$.

  \item[D3] For any $\cF$-measurable, $p$-integrably bounded, and weak$^*$~compact valued correspondence $F$ from $T$ to $X^*$ with $1 < p \leq \infty$, $G_{F}^{\cT}=G_{w^*-\overline{co}F}^{\cT}$.

  \item[D4] Suppose that $G$ is an $\cF$-measurable, $p$-integrably bounded, and weak$^*$ compact valued correspondence from $T$ to $X^*$ with $1< p \leq \infty$. Let $F$ be a weak$^*$ closed valued correspondence from $T \times Y$ to $X^*$ ($Y$ is a metric space). If
      \begin{enumerate}
        \item $F(t,y) \subseteq G(t)$ for $\lambda$-almost all $t\in T$;
        \item for any $y \in Y$, $F(\cdot, y)$ is $\cF$-measurable;
        \item for any $t \in T$, $F(t, \cdot)$ is weak$^*$ upper hemicontinuous;
      \end{enumerate}
      then $H(y)=G_{F_{y}}^{\cT}$ is weak$^*$~upper hemicontinuous in $X^*$.
\end{enumerate}
\end{prop}

\citet{Sun1997} proved that the regularity properties (convexity, compactness, and preservation of upper hemicontinuity) hold for Bochner and Gel$^\prime$fand  integral of correspondences
%\footnote{For integral of Banach valued correspondences, \citet{RY1991} proposed to work with measure spaces with the associated $L^{\infty}$ spaces over any non-null measurable set having strictly larger cardinality than that of the continuum.} 
if the underlying probability space is an atomless Loeb space. These results were generalized by \citet{SY2008} and \citet{Podczeck2008} to those over a saturated probability space. All these papers do not consider the conditional expectation of correspondences. In terms of integral of correspondences, our Propositions~\ref{prop-Bochener}-\ref{prop-G integral} generalize the results in those papers in the sense that the $\sigma$-algebra based on a saturated probability space is nowhere equivalent to any countably-generated sub-$\sigma$-algebra.

In various applications, the widely adopted probability spaces are often not saturated (\textit{e.g.}, based on Polish spaces). Thus, the results depending on saturated probability spaces are not applicable. Instead, we can work with the nowhere equivalence condition, and cover the case of non-saturated probability spaces such as the standard probability spaces.

%%%%%%%%%%%%%%%%%%%%%%%%%%%%%%%%%%%%%%%%%%%%%%%%%%%%%%%%%%%%%%%%%%%%%%%%%%%%%%%%%%%%%%%%%%%%%%%%%%%%%%%%%%%%%%%%%%%%%%%%%%%%%%%%%%%%%%%%%%%%%%%%%%%%%%%%%%%%%%%%%%%%%%%%%%%%%%%%%%%%%%%%%%%%%%%%%%%%%%%%%%%%%%%%%%%%%%%%%%%%%%%%%%%%%%%%

\section{Converse results}\label{sec-converse}

In Sections~\ref{sec-Bochner} and \ref{sec-G integral}, we show that the nowhere equivalence condition is sufficient for several desirable properties of Bochner/Gel$^\prime$fand integral and conditional expectation of correspondences. In this section, we shall prove that  the nowhere equivalence condition is also necessary in the sense that it is implied by each of the properties in Theorems~\ref{thm-ce-bochner} and \ref{thm-ce-G}, and Propositions~\ref{prop-Bochener}, \ref{prop-norm integral} and \ref{prop-G integral}. As a result, the nowhere equivalence condition characterizes those desirable properties.

\begin{thm}\label{thm-converse}
The $\sigma$-algebra $\cT$ is nowhere equivalent to $\cF$ under $\lambda$ if given any sub-$\sigma$-algebra $\cG \subseteq  \cF$, one of the following conditions holds.

\begin{enumerate}
  \item[E1] For any $\cF$-measurable and weakly compact valued (resp. weak$^*$ compact valued) correspondence $F$ from $T$ to $X$ (resp. $X^*$), $CI_{F}^{(\cT,\cG)}$ (resp. $CG_{F}^{(\cT,\cG)}$) is convex.

  \item[E2] For any $\cF$-measurable, $p$-integrably bounded, and weakly (resp. weak$^*$) compact valued correspondence $F$ from $T$ to $X$ (resp. $X^*$) with $1 \le p < \infty$ (resp. $1 < p \le \infty$), $CI_{F}^{(\cT,\cG)}$ (resp. $CG_{F}^{(\cT,\cG)}$) is weakly (resp. weak$^*$) compact in $L_p^\cG(T,X)$ (resp. in the dual space of $L_q^\cG(T,X)$ with $\frac{1}{p} + \frac{1}{q} = 1$).

  \item[E3] For any $\cF$-measurable, $p$-integrably bounded, and weakly (resp. weak$^*$) compact valued correspondence $F$ from $T$ to $X$ (resp. $X^*$) with $1 \le p < \infty$ (resp. $1 < p \le \infty$), $CI_{F}^{(\cT,\cG)}=CI_{\overline{co}F}^{(\cT,\cG)}$ (resp.  $CG_{F}^{(\cT,\cG)} = CG_{w^*-\overline{co}F}^{(\cT,\cG)}$).

  \item[E4] For any weakly (resp. weak$^*$) closed valued correspondence $F$ from $T \times Y$  to $X$ (resp. $X^*$) such that there is an $\cF$-measurable, $p$-integrably bounded, weakly (resp. weak$^*$) compact valued correspondence $G$ from $T$ to $X$ (resp. $X^*$) with $1 \le p < \infty$ (resp. $1 < p \le \infty$),
      \begin{enumerate}
        \item $F(t,y) \subseteq G(t)$ for $\lambda$-almost all $t \in T$;

        \item for any $y \in Y$, $F(\text{\ensuremath{\cdot}},y)$ is $\cF$-measurable;

        \item for any $t \in T$, $F(t,\cdot)$ is weakly (resp. weak$^*$) upper hemicontinuous;
      \end{enumerate}
      $H(y) = CI_{F_y}^{(\cT,\cG)}$ (resp. $CG_{F_y}^{(\cT,\cG)}$) is weakly (resp. weak$^*$) upper hemicontinuous in $L_p^\cG(T,X)$ (resp. in the dual space of $L_q^\cG(T,X)$).

  \item[E5] For any $\cF$-measurable, integrably bounded, and norm compact valued correspondence $F$ from $T$ to $X$, $I_{F}^{\cT}$ is norm compact.

  \item[E6] For any norm closed valued correspondence $F$ from $T\times Y$  to $X$ such that there is an $\cF$-measurable, integrably bounded, and norm compact valued correspondence $G$ from $T$ to $X$,
      \begin{enumerate}
        \item $F(t,y)\subseteq G(t)$ for $\lambda$-almost all $t\in T$;

        \item for any $y\in Y$, $F(\text{\ensuremath{\cdot}},y)$ is $\cF$-measurable;

        \item for any $t\in T$, $F(t,\cdot)$ is norm upper hemicontinuous;
      \end{enumerate}
      $H(y)=I_{F_{y}}^{\cT}$ is norm upper hemicontinuous.
\end{enumerate}
\end{thm}

For the Bochner/Gel$^\prime$fand integral of Banach valued correspondences, the necessity of saturation for those regularity properties was indicated in \citet{SY2008} and \citet{Podczeck2008}. These papers do not distinguish the measurability of a correspondence and its selections, and hence cannot work with the non-saturated probability spaces. As explained in \citet{SY2008}, to show the necessity of saturation, one can construct a particular correspondence with a non-saturated probability space, which is obtained by modifying a counterexample based on the Lebesgue measure space via a measure-preserving mapping. To obtain the necessity of nowhere equivalence, we need to construct a sequence of correspondences satisfying the regularity properties, which is considerably more complicated.

%%%%%%%%%%%%%%%%%%%%%%%%%%%%%%%%%%%%%%%%%%%%%%%%%%%%%%%%%%%%%%%%%%%%%%%%%%%%%%%%%%%%%%%%%%%%%%%%%%%%%%%%%%%%%%%%%%%%%%%%%%%%%%%%%%%%%%%%%%%%%%%%%%%%%%%%%%%%%%%%%%%%%%%%%%%%%%%%%%%%%%%%%%%%%%%%%%%%%%%%%%%%%%%%%%%%%%%%%%%%%%%%%%%%%%%%

\section{Large games}\label{sec-game}

In this section, we apply the results on conditional expectation of Banach valued correspondences to study the existence of pure-strategy Nash equilibria in large games. Similarly as those presented in Sections~\ref{sec-Bochner}, \ref{sec-G integral} and \ref{sec-converse}, we shall state the equilibrium existence results in terms of both Bochner and Gel$^\prime$fand integrable strategy profiles. In particular, the nowhere equivalence condition on the player space is a sufficient and necessary condition for the existence of Nash equilibria in both cases.

A large game can be described as follows.
\begin{enumerate}
\item The player space is an atomless probability space $(T, \cT, \lambda)$. Let $\cF$ be a countably generated sub-$\sigma$-algebra of $\cT$.

\item The set $A$ is the common action space for all players, which is a compact subset of an infinite-dimensional Banach space $X$. Denote $\overline{co}(A)$ as the closed convex hull of $A$, and $\Gamma$ the set of all $\cF$ measurable mapping from $T$ to $\overline{co}(A)$.\footnote{We have not specified the topologies on $A$ and $\Gamma$, which will be explained in the sequel.}

\item The space $\cU$ is the set of continuous real-valued mappings on $A \times \Gamma$ endowed with the supremum norm topology.

\item A game $G$ is an $\cF$-measurable mapping from $T$ to $\mathcal{U}$.
\end{enumerate}

The $\sigma$-algebra $\cF$ here is understood as the one generated by the payoff function $G$ that specifies the players' characteristics. \citet{HSS2017} called $(T, \cF, \lambda)$ the characteristic type space. A strategy profile $g$ is a $\cT$-measurable mapping from the player space $T$ to the action space $A$. The payoff of player~$t \in T$ is $G(t)(g(t), E(g|\cF))$. That is, the payoff of player~$t$ is determined by his own action $g(t)$ as well as $E(g|\cF)$, the societal aggregate given players' characteristics.

%The Borel $\sigma$-algebra on $\mathcal{U}$ is denoted by $\cB(\mathcal{U})$.

\begin{defn}\label{defn-equi}
A Nash equilibrium of a game $G$ is a $\cT$-measurable mapping $g$ from $T$ to $A$ such that for almost all $t\in T$,
$$G(t)(g(t), E(g|\cF)) \ge G(t)(a, E(g|\cF))
$$
for any $a\in A$.%\footnote{We focus on pure-strategy equilibrium.}
\end{defn}

Below, we shall describe both Bochner and Gel$^\prime$fand integrable strategy profiles, and show that the equilibrium existence result can be characterized via the nowhere equivalence condition on the player space. The key point is that there are many players associated with every characteristic; that is, the player space is richer than the characteristic type space for any non-trivial subset of players.

%\begin{defn}\label{defn-atomless}
%If $\cF$ and $\mathcal{G}$ are both $\sigma$-algebras on $T$ and $\mathcal{G}\subseteq\cF$, we say $\cF$ is atomless over $\mathcal{G}$ if for every $D\in\cF$ such that $P(D)>0$, there exists $D_{0}\in\cF$ and $D_{0}\subseteq D$, such that on some set of positive probability,
%\[
%0<P(D_{0}\vert\mathcal{G})<P(D\vert\mathcal{G}).
%\]
%A $\sigma$-algebra is said to be $\aleph_{1}$-atomless if it is atomless over every countably generated $\sigma$-algebra.
%\end{defn}

\textbf{Bochner integrable strategy profile}

\begin{enumerate}
\item The space $X$ is separable, and $X^*$ has the Radon-Nikod\'ym property.

\item The set $A$ is weakly compact, and $\overline{co}(A)$ is the weak closure of the convex hull of $A$.

\item The set $\Gamma$ is viewed as a subset of $L_1^{\cF}(T, X)$ endowed with the weak topology.
\end{enumerate}

\textbf{Gel$^\prime$fand integrable strategy profile}
\begin{enumerate}
\item Let $X = Y^*$, where $Y$ is an infinite-dimensional separable Banach space.

\item The set $A$ is weak$^*$ compact, and $\overline{co}(A)$ is the weak$^*$ closure of the convex hull of $A$.

\item The set $\Gamma$ is viewed as a subset of $L_{\infty}^{\cF}(T, Y^*)$ endowed with the weak$^*$ topology.
\end{enumerate}

\begin{thm} \label{thm-existence}
In each of the cases above, the following statements are equivalent.
\begin{enumerate}
\item The $\sigma$-algebra $\cT$ is nowhere equivalent to $\cF$ under $\lambda$.

\item A Nash equilibrium exists for any large game.
\end{enumerate}
\end{thm}

\begin{rmk}\label{rmk-game-norm}
Theorem~\ref{thm-existence} can be extended to cover the case that the action space $A$ is endowed with the norm topology, where the necessary changes are as follows.
\begin{enumerate}
\item The space $X$ is separable.

\item The set $A$ is norm compact, and $\overline{co}(A)$ is the norm closure of the convex hull of $A$.

\item The space $\mathcal{U}$ is the set of continuous real-valued mappings on $A \times \overline{co}(A)$ endowed with the supremum norm topology. Given a strategy profile $g$, the payoff of player~$t \in T$ is $G(t)(g(t), \int_T g \rmd \lambda)$, where the integral is the Bochner integral. That is, the payoff of player~$t$ is determined by his own action $g(t)$ as well as the societal aggregate $\int_T g \rmd \lambda$.
\end{enumerate}
That is, the externality of the game (the societal aggregate) is the integral of the strategy profile, rather than its conditional expectation. We shall explain how to modify the proof of Theorem~\ref{thm-existence} to cover this case in Remark~\ref{rmk-game-norm proof} of Section~\ref{subsec-game proof}.
\end{rmk}

\citet{KRS1997} presented large games where the societal aggregates are formalized, in one way or another, as Bochner/Gel$^\prime$fand integral. They showed that those games do not possess any Nash equilibrium. \citet{KS1999} worked with a hyperfinite Loeb player space and resolved this nonexistence problem. \citet{SZ2015} proved that the saturation property on the player space is both sufficient and necessary for the existence of Nash equilibrium in large games with infinite-dimensional Banach space as the action space. In all these papers, the societal aggregates are formalized as the integral of the strategy profile.\footnote{\citet{RSY1995} presented a counterexample of large game where the societal aggregate is the distribution induced by the strategy profile. It was then shown in \citet{KS2009} that the saturation property is necessary and sufficient for the existence of Nash equilibrium in such games. \citet{KZ2012} presented a countably generated Lebesgue extension, which is nowhere equivalent to the relevant Borel $\sigma$-algebra, as the player space. Based on this particular Lebesgue extension, they showed that the example in \citet{RSY1995} possesses Nash equilibria. Based on the nowhere equivalence condition, \citet{HSS2017} extended the characterization result of \citet{KS2009}. To prove the necessity result, \citet{HSS2017} constructed a sequence of large games that have the same player space but different action spaces. This result was improved in \citet{HS2018a}, which could work with large games with any specific uncountable action space. In all these papers, the societal aggregate is formalized as the distribution induced by the strategy profile. Note that a player's pure strategy can be viewed as a Dirac measure at that action, and hence the distribution induced by the strategy profile on the action space can be also viewed as the Gel$^\prime$fand integral.} Theorem~\ref{thm-existence} above extends the results in \citet{KS1999} and \citet{SZ2015} in the sense that we can work with non-saturated player spaces and more general societal aggregates (the conditional expectation of the strategy profile). \citet{HS2018b} considered a class of larges games with finite-dimensional action spaces, where the societal aggregate is formalized as the conditional expectation of the strategy profile. Here we consider games with infinite-dimensional action spaces.

%%%%%%%%%%%%%%%%%%%%%%%%%%%%%%%%%%%%%%%%%%%%%%%%%%%%%%%%%%%%%%%%%%%%%%%%%%%%%%%%%%%%%%%%%%%%%%%%%%%%%%%%%%%%%%%%%%%%%%%%%%%%%%%%%%%%%%%%%%%%%%%%%%%%%%%%%%%%%%%%

{\small

\renewcommand{\baselinestretch}{0.7}

\section{Appendix}\label{sec-appendix}

\subsection{Mathematical preliminary}\label{sec-math}

Recall that $(T,\cT,\lambda)$ is an atomless probability space and $\cF$ is a countably generated sub-$\sigma$-algebra of $\cT$. For a correspondence $F$ from $T$ to the Polish space $X$, let
$$D_{F}^{\cT}=\{\lambda \circ f^{-1} \in\cM(X): f\text{ is a }\cT\text{-measurable  selection of }F\},$$
where $\cM(X)$ is the set of Borel probability measures on $X$ endowed with the weak topology, and $\lambda \circ f^{-1}$ is the probability measure on $X$ induced by the mapping $f$. Denote $C_b(X)$ as the set of all bounded continuous function from $X$ to $\bR$. A function $c$ is said to be a Carath\'eodory function if $c(\cdot,x)$ is $\cF$-measurable for each $x\in X$ and $c(t,\cdot)$ is continuous for each $t\in T$.

The following definition introduces the notion of transition probability.
\begin{defn}\label{defn-young measure}
An $\cF$-measurable transition probability from $T$ to $X$ is a mapping $\phi:T\to \cM(X)$ such that $\phi(\cdot,B): t\to \phi(t,B)$ is $\cF$-measurable for every $B\in\cB(X)$.

The set of all $\cF$-measurable transition probabilities from $T$ to $X$ is denoted by $\cR^\cF(X)$, or $\cR^\cF$ when it is clear.
\end{defn}

\begin{defn}\label{defn-convergence rcd}
A sequence $\{\phi_n\}$ in $\cR^\cF$ is said to weakly converge to some $\phi\in\cR^\cF$ if for every bounded Carath\'eodory function $c: T\times X\to \bR$
      $$
      \lim_{n\to\infty}\int_T [\int_X c(t,x) \phi_n(t,\rmd x)]\rmd \lambda(t)= \int_T[\int_X c(t,x) \phi(t,\rmd x)] \rmd\lambda(t).
      $$
The weak topology on $\cR^\cF$ is defined as the weakest topology for which the functional $\phi\to  \int_T[\int_X c(t,x) \phi(t,\rmd x)] \rmd\lambda(t)$  is continuous  for every  bounded Carath\'eodory function $c$.
\end{defn}

If $f$ is a $\cT$-measurable mapping from $T$ to $X$, then $\lambda^{f|\cF}$ denotes the regular conditional distribution (RCD) of $f$ given $\cF$ under $\lambda$. That is, $\lambda^{f|\cF}$ is a mapping from $T\times \cB(X)$ to $[0,1]$ such that
%\footnote{Since $X$ is a Polish space endowed with the Borel $\sigma$-algebra, the RCD $\lambda^{f|\cF}$ always exists; see \citet[Theorem~5.1.9]{Durrett2010}.}
\begin{enumerate}
  \item $\lambda^{f|\cF}(t,\cdot)$ is a probability distribution on $X$ for any $t\in T$;
  \item given any Borel subset $B\subseteq X$, $\lambda^{f|\cF}(\cdot, B)=E[1_B(f)|\cF]$ $\lambda$-a.e.
\end{enumerate}
For a correspondence $F$ from $T$ to $X$, let
$$\cR^{(\cT,\cF)}_F=\{\lambda^{f|\cF}: f \mbox{ is a } \cT\mbox{-measurable seclection of } F\},
$$
which is the set of all RCDs induced by $\cT$-measurable selections of $F$ conditional on $\cF$ under $\lambda$.

The following lemma is from \citet[Theorem~2]{HS2018a}.

\begin{lem}\label{lem-distribution}
If $\cT$ is nowhere equivalent to $\cF$ under $\lambda$, then we have the following properties.
\begin{enumerate}
%  \item For any closed valued $\cF$-measurable correspondence $F$ from $T$ to $X$, $D_{F}^{\cT}$ is convex.
%
%  \item For any compact valued $\cF$-measurable correspondence $F$ from $T$ to $X$, $D_{F}^{\cT}$ is weakly compact.
%
%  \item Let $F$ be a compact valued correspondence from $T$ to X, $Y$ is a metric space, $G$ is a closed valued correspondence from $T\times Y$ to $X$, such that
%     \begin{enumerate}
%        \item[a] $\forall(t,y)\in T\times Y$, $G(t,y)\subseteq F(t)$;
%        \item[b] $\forall y\in Y$, $G(\text{\ensuremath{\cdot}},y)$ (denoted as $G_y$) is $\cF$-measurable from $T$ to $X$;
%        \item[c] $\forall t\in T$, $G(t,\cdot)$ (denoted as $G_t$) is upper hemicontinuous from $Y$ to $X$;
%      \end{enumerate}
%      Then $H(y)=D_{G_{y}}^{\cT}$ is weakly upper hemicontinuous from $Y$ to $\cM(X)$.

  \item For any closed valued $\cF$-measurable correspondence $F$ from $T$ to $X$, $\cR^{(\cT,\cG)}_F$ is convex.

  \item For any compact valued $\cF$-measurable correspondence $F$ from $T$ to $X$, $\cR^{(\cT,\cG)}_F$ is weakly compact.

  \item Let $F$ be a compact valued $\cF$-measurable correspondence from $T$ to $X$, and $G$ a closed valued correspondence from $T\times Y$ to $X$ ($Y$ is a metric space) such that
      \begin{enumerate}
        \item[a] for any $(t,y)\in T\times Y$, $G(t,y)\subseteq F(t)$;
        \item[b] for any $y\in Y$, $G(\text{\ensuremath{\cdot}},y)$ (denoted as $G_y$) is $\cF$-measurable from $T$ to $X$;
        \item[c] for any $t\in T$, $G(t,\cdot)$ (denoted as $G_t$) is upper hemicontinuous from $Y$ to $X$.
      \end{enumerate}
      The correspondence $H(y)=\cR_{G_{y}}^{(\cT,\cG)}$ is upper hemicontinuous from $Y$ to $\cR^{\cG}$.
  \end{enumerate}
\end{lem}

\begin{rmk} \label{rmk-distribution}
Note that Lemma~\ref{lem-distribution} above is a corollary of Theorem~\ref{thm-ce-G}. Every Polish space admits a compatible totally bounded metric. Let $\hat{X}$ be the compactification of $X$ under this compatible totally bounded metric. Then $C_b(\hat{X})$, the space of bounded continuous functions on $\hat{X}$, is a Polish space. The dual space of $C_b(\hat{X})$ is the space of finite measures on $\hat{X}$.

The correspondence $F$ in Lemma~\ref{lem-distribution} can be also viewed as the collection of all Dirac measures at its values. That is, one can equivalently work with the correspondence $\tilde{F}(t) = \{\delta_x \colon x \in F(t) \}$. Then $\cR^{(\cT,\cG)}_F$, the set of all RCDs induced by $\cT$-measurable selections of $F$ conditional on $\cG$, is the Gel$^\prime$fand conditional expectation of $\tilde{F}(t)$. Even though $\tilde{F}$ is a correspondence from $T$ to the space of finite measures on $\hat{X}$, it only takes values from the space of finite measures on $X$. Lemma~\ref{lem-distribution} now follows from (D1), (D2) and (D4) of Theorem~\ref{thm-ce-G}.
 \end{rmk}

Below, we introduce a technical condition, which will be used extensively for deriving the necessity results in Theorems~\ref{thm-converse} and \ref{thm-existence}. Lemma~\ref{lem-neq} is part of \citet[Lemma~1]{HS2018a}.

\begin{defn}\label{defn-others}
The $\sigma$-algebra $\cF$ admits an \textbf{asymptotic independent supplement} in $\cT$ if for some strictly increasing sequence $\{n_k\}_{k=1}^\infty$ of positive integers and for each $k \ge 1$, there exists a $\cT$-measurable partition $\{E_1,E_2,\ldots,E_{n_k}\}$ of $T$ with $\lambda(E_j) = \frac{1}{n_k}$ and $E_j$ being independent of $\cF$ for $j=1,2,\ldots, n_k$.
\end{defn}

\begin{lem}\label{lem-neq}
The following statements are equivalent.
\begin{enumerate}[(i)]
\item $\cT$ is nowhere equivalent to $\cF$.
\item $\cF$ admits an asymptotic independent supplement in $\cT$.
\end{enumerate}
\end{lem}

%%%%%%%%%%%%%%%%%%%%%%%%%%%%%%%%%%%%%%%%%%%%%%%%%%%%%%%%%%%%%%%%%%%%%%%%%%%%%%%%%%%%%%%%%%%%%%%%%%%%%%%%%%%%%%%%%%%%%%%%%%%%%%%%%%%%%%%%%%%%%%%%%%%%%%%%%%%%%%%%

\subsection{Proofs of the results in Section~\ref{sec-Bochner}}\label{proof-Bochner}

%\subsubsection{Proofs of Subsection~\ref{subsec-regularity}}

To prove Theorem~\ref{thm-ce-bochner}, we first provide several useful lemmas.

Recall that $X$ is a separable Banach space. The unit ball of $X^*$ is metrizable under the weak$^*$  topology (see \citet[Theorem 6.30]{AB2006}), and hence separable. There is a countable dense set $\{x^*_m\}_{m\in\bN}$ in the unit ball of $X^*$. Fix a sub-$\sigma$-algebra $\cG$ of $\cF$ and an $\cF$-measurable correspondence $F$ from $T$ to $X$. Suppose that $f$ is a $\cT$-measurable Bochner integrable selection of $F$.

Let $\Phi$ be the set of all mappings $\varphi$ from $(T,\cG,\lambda)$ to $X^*$ such that $\varphi(t)\in \{x_m^*\}_{m\in \bN}$ for $\lambda$-almost all $t\in T$. The following lemma shows that to distinguish distinct elements in $L^{\cG}_p(T,X)$, instead of working with the dual space $L_q^\cG(T,X^*)$ with $\frac{1}{p}+\frac{1}{q}=1$, it is sufficient to consider the set $\Phi$.

\begin{lem}\label{lem-separate}
The set $\Phi$ separates points in $L_p^\cG(T,X)$ for $1 \le p < \infty$.
\end{lem}

\begin{proof}
Recall that $X^*$ has the Radon-Nikod\'ym property.  Since $X$ is separable and $\cG$ is countably generated, $L^{\cG}_p(T,X)$ is also a separable Banach space for $1 \leq p <\infty$. By Theorem 1 in \citet[p.98]{DU1977}, $L^{\cG}_{q}(T, X^*)$ consists of continuous linear functionals of $L_{p}^\cG(T,X)$  for $1 \leq p < \infty$, where $\frac{1}{p}+\frac{1}{q}=1$. Given two different mappings $f_1, f_2\in L_p^\cG(T,X)$, there exists a constant $\epsilon > 0$ and a mapping $g\in L_q^\cG(T,X^*)$ such that
$$\int_T g(t, f_1(t)) \lambda(\rmd t) - \int_T g(t, f_2(t)) \lambda(\rmd t) > \epsilon.\footnote{Given a mapping $g$ from $T$ to $X^*$ and $x \in X$, we slightly abuse the notation by using $g(t, x)$ to denote $g(t)(x)$ when there is no confusion.}
$$
Without loss of generality, we assume that $g\in L_\infty^\cG(T,X^*)$, and is essentially bounded by $1$.

For each $m\in \bN$, let
$$D_m = \left\{ t\in T\colon \left|x_m^*(f_i(t)) -g(t, f_i(t))\right| < \frac{\epsilon}{3}, i=1,2 \right\}.
$$
Since $\{x_m^*\}_{m\in \bN}$ is dense under the weak$^*$ topology in the unit ball of $X^*$, $\lambda(\cup_{m\in \bN} D_m) = 1$. Define a mapping $\varphi \in \Phi$ such that $\varphi(t) = x_m^*$ for $t\in D_m\setminus \left( \cup_{0\le k \le m-1} D_k \right)$, where $D_0 =\emptyset$. Then
$$ \int_T  \varphi f_1 \rmd \lambda - \int_T  \varphi f_2 \rmd \lambda  \ge \int_T  (g f_1 - \frac{\epsilon}{3}) \rmd \lambda - \int_T  (g f_2 + \frac{\epsilon}{3}) \rmd \lambda > \frac{\epsilon}{3}.
$$
This completes the proof.
\end{proof}

The following result has been shown in \citet[proof of Theorem 3]{Sun1997}, which will be useful for proving our Theorem~\ref{thm-ce-bochner}. We summarize the result as Lemma~\ref{lem-sun1997} and provide the proof for completeness.

\begin{lem}\label{lem-sun1997}
Let
$$\rho_w(x,y) = \sum^\infty_{m=1} {1\over{2^m}} |x^*_m (x - y)|
$$ for any $x,y \in X$. The Borel $\sigma$-algebras on $X$ generated by the norm topology, the weak topology and the topology induced by the metric $\rho_w$ are the same.
\end{lem}

\begin{proof}
Note that on $X$, the topology induced by the metric $\rho_w$ is weaker than the weak topology (see Theorem~6.3 in \citet[p.434]{DS1958}). This implies that the identity mapping $I_X$ is an injective continuous mapping from $X$ (endowed with the norm topology) to $X$ (endowed with the topology induced by $\rho_w$).

Let $(Z,\rho_w)$ be the completion of the metric space $(X,\rho_w)$. Then $I_X$ is a Borel measurable injective mapping from the Polish space $X$ with the norm topology to the Polish space $(Z,\rho_w)$. By Corollary~3.3 in \citet[p. 22]{Parthasarathy1967}, $X$ is a Borel subset of $(Z, \rho_w)$ and $I_X$ is a Borel isomorphism from the space $(X, \|\cdot\|)$ to the space $(X,\rho_w)$. Thus, the Borel $\sigma$-algebras on $X$ generated by the norm topology, the weak topology and the topology induced by $\rho_w$ are all the same. In addition, since the identity mapping $I_X$ is a continuous mapping from $X$ (endowed with the weak/norm topology) to the metric space $(Z,\rho_w)$, any weakly/norm compact set in $X$ is also compact in $(Z,\rho_w)$.
\end{proof}

Next, we present an equivalent version of the conditional expectation $E(f|\cG)$ of $f$ conditional on $\cG$, which is easier to work with. Given a Bochner integrable mapping $f$ from $(T, \cT, \lambda)$ to $X$ (endowed with the weak topology), by Lemma~\ref{lem-sun1997} above, one can also view $f$ as a measurable mapping from $(T, \cT, \lambda)$ to the Polish space $(Z, \rho_w)$.

\begin{lem}\label{lem-Bochner ce}
For a Bochner integrable mapping $f$ from $(T, \cT, \lambda)$ to $X$,
$$E(f|\cG)(t) = \int_X I_{X}(x) \lambda^{f|\cG}(t,\rmd x)
$$
for $\lambda$-almost all $t\in T$, where (1) $I_X(\cdot)$ is the identity mapping from $X$ (endowed with the topology induced by the metric $\rho_w$) to $X$ (endowed with the weak topology), and (2) $\lambda^{f|\cG}$ is the regular conditional distribution given $\cG$ generated by viewing $f$ as a mapping from $(T, \cT, \lambda)$ to the Polish space $(Z, \rho_w)$.
\end{lem}

\begin{proof}
By Theorem~4.38 in \citet{AB2006}, the mapping $I_X(\cdot)$ can be approximated by a sequence of simple functions $\{\psi_k\}_{k \ge 1}$ pointwise. In particular,  the sequence $\{\psi_k\}_{k \ge 1}$ can be chosen such that for each $k \ge 1$,
\begin{enumerate}
\item $\psi_k(x) = \sum_{1 \le j \le n_k}1_{C_{kj}}(x) x_{kj}$, where $n_k \ge 1$ is some positive integer, $\{x_{kj}\}_{1 \le j \le n_k}$ are distinct points in $X$, and $\{ C_{kj} \}_{1 \le j \le n_k}$ is a partition of $X$;

\item $\|\psi_k(x) \| \le 2 \|x\|$ for each $x \in X$.
\end{enumerate}

By Lemma~\ref{lem-sun1997}, $f$ is also measurable if one views it as a mapping from $(T, \cT, \lambda)$ to $(Z, \rho_w)$. In particular, it takes values from $(X, \rho_w)$. By the dominated convergence theorem (see Theorem~11.46 in \citet{AB2006}), for any $t$, $\int_X \psi_k(x) \lambda^{f|\cG}(t,\rmd x) \to \int_X I_{X}(x) \lambda^{f|\cG}(t,\rmd x)$ as $k \to \infty$. It further implies that for any $\phi \in \Phi$,
$$\int_{T} \phi \left(t, \int_X \psi_k(x) \lambda^{f|\cG}(t,\rmd x) \right) \lambda(\rmd t)  \to \int_{T} \phi \left(t, \int_X I_X(x) \lambda^{f|\cG}(t,\rmd x) \right) \lambda(\rmd t)
$$
as $k \to \infty$. In addition,
\begin{align*}
\int_{T} \phi \left(t, \int_X \psi_k(x) \lambda^{f|\cG}(t,\rmd x) \right) \lambda(\rmd t)
& = \int_{T} \phi \left(t, \int_X \sum_{1 \le j \le n_k} 1_{C_{kj}}(x) x_{kj} \lambda^{f|\cG}(t,\rmd x) \right) \lambda(\rmd t) \\
& = \sum_{1 \le j \le n_k} \int_{T} \phi (t, x_{kj}) \lambda^{f|\cG}(t,C_{kj}) \lambda(\rmd t) \\
& = \sum_{1 \le j \le n_k} \int_{T} \phi (t, x_{kj}) E(1_{C_{kj}}(f)|\cG)(t) \lambda(\rmd t) \\
& = \sum_{1 \le j \le n_k} \int_{T} E\left( \phi (t, x_{kj}) \cdot 1_{C_{kj}}(f)|\cG \right) (t) \lambda(\rmd t) \\
& = \sum_{1 \le j \le n_k} \int_{T} \phi (t, x_{kj}) 1_{C_{kj}}(f(t)) \lambda(\rmd t) \\
& = \int_{T} \phi (t, \psi_k(f(t))) \lambda(\rmd t) \\
& \to \int_{T} \phi (t, f(t)) \lambda(\rmd t).
\end{align*}
The third equality holds due to the definition of regular conditional distribution. The fourth equality is true since $\phi$ is $\cG$-measurable. The last line again follows from the dominated convergence theorem. The other equalities are just simple algebras.

As a result, we have that for any $\phi \in \Phi$,
\begin{align*}
\int_{T} \phi \left(t, \int_X I_X(x) \lambda^{f|\cG}(t,\rmd x) \right) \lambda(\rmd t)
& = \int_{T} \phi (t, f(t)) \lambda(\rmd t) \\
& = \int_{T} \phi (t, E(f|\cG)(t)) \lambda(\rmd t).
\end{align*}
By Lemma~\ref{lem-separate}, the set $\Phi$ separates points in $L_1^\cG(T,X)$, which implies that
$$E(f|\cG)(t) = \int_X I_{X}(x) \lambda^{f|\cG}(t,\rmd x)
$$
for $\lambda$-almost all $t\in T$.
\end{proof}

Now we are ready to prove Theorem~\ref{thm-ce-bochner}.

\begin{proof}[Proof of Theorem~\ref{thm-ce-bochner}]

\

A1. Recall that $(Z,\rho_w)$ is the completion of the metric space $(X,\rho_w)$. Then $F$ can be viewed as an $\cF$-measurable and compact valued correspondence from $T$ to $Z$, taking values from $X$. Fix two $\cT$-measurable Bochner integrable selections $f_1$, $f_2$ and $0\leq \alpha \leq 1$. By Lemma~\ref{lem-distribution}~(1), $\cR^{(\cT,\cG)}_{F}$ is convex. Then there exists a $\cT$-measurable selection $f$ of $F$ such that $\lambda^{f|\cG} = \alpha \lambda^{f_1|\cG} + (1-\alpha) \lambda^{f_2|\cG}$. By Lemma~\ref{lem-Bochner ce}, for $\lambda$-almost all $t\in T$,
\begin{align*}
\alpha E(f_1|\cG)(t) + (1-\alpha) E(f_2|\cG)(t)
& = \int_X I_X(x) \left(\alpha \lambda^{f_1|\cG}(t, \rmd x) + (1-\alpha) \lambda^{f_2|\cG}(t, \rmd x) \right) \\
& = \int_X I_X(x) \lambda^{f|\cG}(t,\rmd x) \\
& = E(f|\cG)(t),
\end{align*}
which implies that $CI_{F}^{(\cT,\cG)}$ is convex.

%Suppose that the correspondence $F$ is $\cF$-measurable from $T$ to $X$ and norm closed valued. Fix two $\cT$-measurable Bochner integrable selections $f_1$, $f_2$ and $0\leq \alpha \leq 1$. By Lemma~\ref{lem-distribution}~(1), $\cR^{(\cT,\cG)}_{F}$ is convex. Then there exists a $\cT$-measurable selection $f$ of $F$ such that $\lambda^{f|\cG} = \alpha \lambda^{f_1|\cG} + (1-\alpha) \lambda^{f_2|\cG}$. By Lemma~\ref{lem-Bochner ce}, for $\lambda$-almost all $t\in T$,
%\begin{align*}
%\alpha E(f_1|\cG)(t) + (1-\alpha) E(f_2|\cG)(t)
%& = \int_X I_X(x) \left(\alpha \lambda^{f_1|\cG}(t, \rmd x) + (1-\alpha) \lambda^{f_2|\cG}(t, \rmd x) \right) \\
%& = \int_X I_X(x) \lambda^{f|\cG}(t,\rmd x) \\
%& = E(f|\cG)(t),
%\end{align*}
%which implies that $CI_{F}^{(\cT,\cG)}$ is convex.

%Alternatively, suppose that the correspondence $F$ is $\cF$-measurable from $T$ to $X$ and weakly compact valued. Recall that $(Z,\rho_w)$ is the completion of the metric space $(X,\rho_w)$, and the identity mapping $I_X$ is an injective continuous mapping from $X$ (endowed with the weak topology) to $X$ (endowed with the topology induced by $\rho_w$). By Lemma~\ref{lem-sun1997}, the Borel $\sigma$-algebras on $X$ generated by the weak topology and the topology induced by the metric $\rho_w$ are the same. Thus, $F$ can be viewed as an $\cF$-measurable and compact valued correspondence from $T$ to $Z$, taking values from $X$. One can then repeat the argument in the previous paragraph and complete the proof.

\

A2. By Eberlein-Smulian Theorem (see Theorem~6.34 in \citet{AB2006}), it is sufficient to show that $CI_{F}^{(\cT,\cG)}$ is weakly sequentially compact in $L_p^\cG(T,X)$. That is, given an arbitrary sequence of $\cT$-measurable Bochner integrable selections $\{f_n\}_{n\in \bN}$ of $F$ with $g_n = E(f_n|\cG)$ for each $n\in \bN$, there is a subsequence of $\{g_n\}_{n\in \bN}$ that weakly converges in $L_p^\cG(T,X)$ to some point in $CI_{F}^{(\cT,\cG)}$.

For the case $p=1$, since the sequence $\{f_n\}_{n\in \bN}$ is integrably bounded and $F$ is weakly compact valued, it is relatively weakly compact in $L_1^\cT(T,X)$ (see Corollary~2.6 in \citet{DRS1993}). Without loss of generality, we assume that $f_n$ weakly converges to some $f_0 \in L_1^\cT(T,X)$. Then $g_n$  weakly converges to $E(f_0|\cG)$ in $L_1^\cG(T,X)$.

For the case $1 < p < \infty$, since the $p$-integrably boundedness implies the integrably boundedness, we have that the sequence $\{f_n\}_{n\in \bN}$ is relatively weakly compact in $L_1^\cT(T,X)$. Combining Theorem~2.1 and Corollary~3.4 in \citet{DRS1993}, we have that the sequence $\{f_n\}_{n \in \bN}$ is relatively weakly compact in $L_p^\cT(T,X)$.  Without loss of generality, we assume that $f_n$ weakly converges to some $f_0 \in L_p^\cT(T,X)$. Then $g_n$  weakly converges to $E(f_0|\cG)$ in $L_p^\cG(T,X)$.

%Let $F_1(t) = \overline{\{f_n(t)\}_{n \in \bN}}$ for each $t \in T$, where $\overline{\{f_n(t)\}_{n \in \bN}}$ is the weak closure of the set $\{f_n(t)\}_{n \in \bN}$ in $X$.

Recall that $(Z,\rho_w)$ is the completion of $(X,\rho_w)$. We can view $F$ as a correspondence from $(T, \cF,\lambda)$ to the Polish space $(Z,\rho_w)$. Since $F$ is compact valued, $\cR^{(\cT,\cG)}_{F}$ is compact by Lemma~\ref{lem-distribution}~(2). Thus, there is a $\cT$-measurable Bochner integrable selection $f$ of $F$ and a subsequence of $\{f_{n}\}_{n\in \bN}$, say itself,  such that $\lambda^{f_n|\cG}$ weakly converges to $\lambda^{f|\cG}$ as $n\to \infty$.

%Note that $f$ is also a $\cT$-measurable Bochner integrable selection of $F$.

For each $m\in \bN$, $x_m^*$ is uniformly continuous on $(X,\rho_w)$, and hence has a unique continuous extension $\tilde{x_m^*}$ on $(Z,\rho_w)$. Fix some $\varphi$ from $(T,\cG,\lambda)$ to $X^*$ such that $\varphi(t)\in \{x_m^*\}_{m\in \bN}$ for $\lambda$-almost all $t\in T$. Denote $E_m = \{t\in T \colon \varphi(t) = x_m^*\}$ for $m \ge 1$. Let $\varphi_m (t) = \mathbf{1}_{E_m}(t) \cdot x_m^*$ for $m \ge 1$. Then $\varphi(t) = \sum_{m \in \bN} \varphi_m (t)$ for $\lambda$-almost all $t \in T$, and $\varphi_m$ is $\cG$-measurable since $E_m \in \cG$ for each $m \in \bN$.
%and $E_0 = T \setminus \cup_{m \ge 1} E_m$
We can define a mapping $\phi$ from $(T,\cG,\lambda)$ to the space of continuous functions on $(Z,\rho_w)$ by letting $\phi(t) = \tilde{x_k^*}$ if $\varphi(t) = x_k^*$ for some $k\in \bN$, and $0$ otherwise. For each $m\in \bN$, $\varphi_m$ can be extended to $\phi_m = \mathbf{1}_{E_m}(t) \cdot \tilde{x_m^*}$ similarly.

Fix $m\in \bN$. Denote $\psi_m^n(t) = \phi_m(t,f_n(t))$ and $\psi_m(t) = \phi_m(t,f(t))$.  Given any bounded continuous function $c$ on $\bR$,
\begin{align*}
\int_\bR c(a) \lambda\circ \left(\psi_m^n\right)^{-1} (\rmd a)
& = \int_T c(\phi_m(t,f_n(t))) \rmd\lambda(t) = \int_T c(\mathbf{1}_{E_m}(t) \cdot \tilde{x_m^*}(f_n(t))) \rmd\lambda(t)   \\
& = \lambda(E_m^c)\cdot c(0) + \int_{E_m} c(\tilde{x_m^*}(f_n(t))) \rmd\lambda(t) \\
& = \lambda(E_m^c)\cdot c(0) + \int_{E_m} E\left( c(\tilde{x_m^*}(f_n(t))) |\cG \right) \rmd\lambda(t) \\
& = \lambda(E_m^c)\cdot c(0) + \int_{E_m} \int_Z c(\tilde{x_m^*}(x)) \lambda^{f_n|\cG}(t, \rmd x) \rmd\lambda(t)  \\
& \to \lambda(E_m^c)\cdot c(0) + \int_{E_m} \int_Z c(\tilde{x_m^*}(x)) \lambda^{f|\cG}(t, \rmd x) \rmd\lambda(t) \\
& = \int_\bR c(a) \lambda\circ \left(\psi_m\right)^{-1} (\rmd a)
\end{align*}
Therefore, $\lambda \circ \left(\psi_m^n\right)^{-1}$ weakly converges to $\lambda \circ \psi_m^{-1}$ as $n\to \infty$.

Since $F$ is integrably bounded, the sequence $\{\psi_m^n\}$ is uniformly integrable. By Theorem 5.4 of \citet[p.32]{Billingsley1968}, we have that
$$\lim_{n\to\infty} \int_T \psi_m^n \rmd\lambda  = \int_T \psi_m \rmd\lambda.
$$
By the dominated convergence theorem, we can take summation for all $m\in \bN$ and conclude that
$$\lim_{n\to\infty} \int_T \psi^n \rmd\lambda  = \int_T \psi \rmd\lambda,
$$
where $\psi^n(t) = \phi(t,f_n(t))$ and $\psi(t) = \phi(t,f(t))$. We have that
\begin{align*}
\int_T \psi(t) \rmd\lambda(t)
& = \int_T \phi(t,f(t)) \rmd\lambda(t) = \int_T \varphi(t,f(t)) \rmd\lambda(t) = \sum_{m=1}^\infty \int_T \varphi_m(t,f(t)) \rmd\lambda(t) \\
& = \sum_{m=1}^\infty \int_T \mathbf{1}_{E_m}(t) \cdot x_m^*(f(t)) \rmd\lambda(t) = \sum_{m=1}^\infty \int_T \mathbf{1}_{E_m}(t) \cdot x_m^*(E(f|\cG)(t)) \rmd\lambda(t) \\
& = \int_T \varphi(t,E(f|\cG)(t)) \rmd\lambda(t) = \int_T \varphi(t,g(t)) \rmd\lambda(t),
\end{align*}
where $g = E(f|\cG)$. Similarly, we could show that for each $n\in \bN$,
$$\int_T \psi^n(t) \rmd\lambda(t)  = \int_T \varphi(t,g_n(t)) \rmd\lambda(t).
$$
Thus,
$$\lim_{n\to\infty} \int_T \varphi(t,g_n(t)) \rmd\lambda(t) =  \int_T \varphi(t,g(t)) \rmd\lambda(t).$$
By Lemma~\ref{lem-separate}, the set of all such $\varphi$ separates points in $L_p^\cG(T,X)$. Then we have that $E(f_0|\cG) = g$, which implies that $CI_{F}^{(\cT,\cG)}$ is weakly compact.

\

A3. We shall first show that $\overline{CI_{\overline{co}F}^{(\cT,\cG)}} = \overline{co}\,CI_{F}^{(\cT,\cG)}$, where
$\overline{CI_{\overline{co}F}^{(\cT,\cG)}}$ (resp. $\overline{co}\,CI_{F}^{(\cT,\cG)}$) is the norm closure of $CI_{\overline{co}F}^{(\cT,\cG)}$ (resp. the convex hull of $CI_{F}^{(\cT,\cG)}$) in $L_p^\cG(T,X)$. It is obvious that $\overline{co}\,CI_{F}^{(\cT,\cG)} \subseteq \overline{CI_{\overline{co}F}^{(\cT,\cG)}}$, we only need to show the other direction.

Let $f_0$ be a $\cT$-measurable Bochner integrable selection of $\overline{co} F$ and $g_0 = E(f_0|\cG)$. For any $\varphi\in L^\cG_q(T,X^*)$ with $\frac{1}{p} + \frac{1}{q} = 1$,
\begin{align*}
\varphi(g_0)
& = \int_T \varphi(t,g_0(t)) \rmd\lambda(t) = \int_T \varphi(t, E(f_0|\cG)(t)) \rmd\lambda(t) \\
& = \int_T E(\varphi(t, f_0(t))|\cG) \rmd\lambda(t) = \int_T \varphi(t,f_0(t)) \rmd\lambda(t) \\
& \leq  \int_T \sup_{x\in \overline{co}F(t)}\varphi(t,x) \rmd\lambda(t) = \int_T \sup_{x\in F(t)}\varphi(t,x) \rmd\lambda(t).
\end{align*}
The third equality holds since $\varphi$ is $\cG$-measurable. The inequality is true since $f_0$ is a selection of
$\overline{co} F$. The last equality holds because $\varphi(t)$ is linear in $X$ for each $t \in T$.

Fix $\epsilon >0$ and $w\in L_1^\cT(T,\bR_+)$. Let
$$\Phi(t) = [\sup_{x\in F(t)}\varphi(t,x) - \epsilon w(t), \sup_{x\in F(t)} \varphi(t,x)].
$$
By Theorem 8.2.9 of \citet{AF1990}, there exists a $\cT$-measurable selection $f$ of $F$ such that for $\lambda$-almost all $t\in T$,
$$\varphi(t,f(t)) \geq \sup_{x\in F(t)}\varphi(t,x) - \epsilon w(t).$$
Then we have
\begin{align*}
\int_T \sup_{x\in F(t)}\varphi(t,x) \rmd\lambda
& \le \int_T [\varphi(t,f(t)) + \epsilon w(t)] \rmd\lambda(t) \\
& \le \int_T \varphi(t,f(t)) \rmd\lambda(t) + \epsilon \int_T |w(t)| \rmd \lambda(t) \\
& = \int_T \varphi(t, E(f|\cG)(t)) \rmd\lambda(t) + \epsilon \int_T |w(t)| \rmd \lambda(t) \\
& \le \sup_{g\in CI_{F}^{(\cT,\cG)}} \varphi(g) + \epsilon \int_T |w(t)| \rmd\lambda(t);
\end{align*}
that is, $\varphi(g_0)  \leq \sup_{g\in CI_{F}^{(\cT,\cG)}} \varphi(g)$. Since $\varphi$ is arbitrary, we have that $\varphi(g_0)  \geq \inf_{g\in \overline{co}CI_{F}^{(\cT,\cG)}} \varphi(g)$. Note that $\overline{co}CI_{F}^{(\cT,\cG)}$ is a closed convex set. By the separation theorem, $g_0 \in \overline{co}CI_{F}^{(\cT,\cG)}$, which further implies that $\overline{CI_{\overline{co}F}^{(\cT,\cG)}}  \subseteq  \overline{co}\,CI_{F}^{(\cT,\cG)}$. That is, $\overline{CI_{\overline{co}F}^{(\cT,\cG)}} = \overline{co}\,CI_{F}^{(\cT,\cG)}$.

By (A1) and (A2), $CI_{F}^{(\cT,\cG)}$ and $CI_{\overline{co}F}^{(\cT,\cG)}$ are both convex and weakly compact. By Mazur's Theorem (see \citet[p.~292]{RF2010}), a convex and weakly closed set is also norm closed. Thus, $\overline{co}\,CI_{F}^{(\cT,\cG)} = CI_{F}^{(\cT,\cG)}$ and $\overline{CI_{\overline{co}F}^{(\cT,\cG)}} = CI_{\overline{co}F}^{(\cT,\cG)}$, and hence $CI_{\overline{co}F}^{(\cT,\cG)}  = CI_{F}^{(\cT,\cG)}$.

\

A4. By (A2), $CI_G^{(\cT,\cG)}$ is weakly compact, and hence $CI_{F_y}^{(\cT,\cG)}$ is a weakly compact valued correspondence from $Y$ to the weakly compact set $CI_G^{(\cT,\cG)}$. As shown in the proof of Lemma~\ref{lem-sun1997}, the weak topology on a weakly compact set in a separable Banach space is metrizable. Thus, $CI_{F_y}^{(\cT,\cG)}$ can be viewed as a compact valued correspondence from $Y$ to a compact metric space.

Let $\{y_n\}_{n=0}^\infty$ be a sequence in $Y$ ($y_n \neq y_0$ for any $n \ge 1$), and $\{f_n\}_{n \ge 1}$ a sequence of $\cT$-measurable and Bochner integrable selection of $F_{y_n}$ for $n\ge 1$. Suppose that $g_n$ weakly converges to some $g_0$ in $L_p^\cG(T,X)$ with $g_n = E(f_n|\cG)$ for $n \ge 1$, and $\lim_{n\to\infty} y_n = y_0$. It suffices to show that $g_0 \in CI_{F_{y_0}}^{(\cT,\cG)}$ (see Theorem~17.11 in \citet{AB2006}).

Due to the assumption above, $g_0 \in w-\limsup_{n\to\infty} CI_{F_{y_n}}^{(\cT,\cG)}$.\footnote{For a sequence of sets $\{A_n\}_{n \ge 1}$, (1) $\limsup_{n \to \infty} A_n$ is the topological limit supermum of $\{A_n\}$, which collects all the points $x$ such that every neighborhood of $x$ intersects infinitely many $A_n$; (2) $\liminf_{n\to\infty} A_n$ is the topological limit infimum of $\{A_n\}$, which collects all the points $x$ such that every neighborhood of $x$ intersects all the $A_k$ for sufficiently large $k$. The notation ``$w-$'' means that it is under the weak topology.} Since the correspondence $F(t, \cdot)$ is weakly upper hemicontinuous for each $t\in T$,  $w-\limsup_{n\to\infty} F(t, y_n) \subseteq F(t, y_0)$. Recall that the Polish space $(Z,\rho_w)$ is the completion of $(X,\rho_w)$. As mentioned above, the weak topology on a weakly compact set in a separable Banach space is metrizable under the metric $\rho_w$. Thus, the topological limit supermum of the sequence $\{F_{y_n}(t)\}_{n=1}^\infty$ (included in the compact set $G(t)$) in $(Z, \rho_w)$ is the same as $w-\limsup_{n\to\infty} F_{y_n}(t)$. By Lemma~\ref{lem-distribution}~(3), there is a subsequence of $\{f_{n}\}^\infty_{n=1}$, say  itself, and a $\cT$-measurable selection $f$ of $F_{y_0}$ such that $\lambda^{f_{n}|\cG}$ weakly converges to $\lambda^{f|\cG}$.

As in the proof of (A2), it is straightforward to show that $\lim_{n\to \infty}\varphi(g_n) = \varphi(g)$ for each $\varphi \in \Phi$, where $g = E(f|\cG)$. By Lemma~\ref{lem-separate}, the set $\Phi$ separates points in $L_p^\cG(T,X)$, which implies that $g_0 = g \in CI_{F_{y_0}}^{(\cT,\cG)}$.
\end{proof}

The following result is a straightforward corollary of Theorem~\ref{thm-ce-bochner}~(A4).

\begin{coro}\label{coro-CI weak upper}
Suppose that $G$ is an $\cF$-measurable, $p$-integrably bounded and weakly compact valued correspondence from $T$ to $X$ for $1 \le p < \infty$. Let $\{F_n\}^\infty_{n=1}$ be a sequence of $\cF$-measurable correspondences from $T$ to $X$ such that $F_n(t)\subseteq G(t)$ for $\lambda$-almost all $t\in T$ and $n\geq 1$. Let $F_0 = w-\limsup_{n \to \infty} F_n$. Then $w-\limsup_{n\to\infty} CI_{F_n}^{(\cT,\cG)} \subseteq CI_{F_0}^{(\cT,\cG)}$.
\end{coro}

In Theorem~\ref{thm-ce-bochner}~(A2) and (A4), we study the properties on the compactness and preservation of upper hemicontinuity of conditional expectation of correspondences via the condition of nowhere equivalence. Below, we show that under some regularity conditions, the properties on the compactness and preservation of upper/lower hemicontinuity for conditional expectation of correspondences hold without imposing the nowhere equivalence condition. Note that in Proposition~\ref{prop-CI-hemicontinuity}~(F3), we do not require the correspondence $F$ to be convex valued for the preservation of lower hemicontinuity.

\begin{prop}\label{prop-CI-hemicontinuity}
For any sub-$\sigma$-algebra $\cG$ of $\cF$, we have the following properties.
\begin{enumerate}
\item[F1] For any $\cF$-measurable, $p$-integrably bounded, convex and weakly compact valued correspondence $F$ from $T$ to $X$ with $1\leq p < \infty$, $CI_{F}^{(\cT,\cG)}$ is convex and weakly compact in $L_p^\cG(T,X)$.

\item[F2] Suppose that $G$ is an $\cF$-measurable, $p$-integrably bounded and weakly compact valued correspondence from $T$ to $X$ with $1 \le p < \infty$. Let $F \colon T\times Y\to X$ be a convex and weakly closed valued correspondence ($Y$ is a metric space). If
      \begin{enumerate}
        \item $F(t,y)\subseteq G(t)$ for $\lambda$-almost all $t\in T$;
        \item for any $y\in Y$, $F(\text{\ensuremath{\cdot}},y)$ is $\cF$-measurable;
        \item for any $t\in T$, $F(t,\cdot)$ is weakly upper hemicontinuous;
      \end{enumerate}
      then $H(y)= CI_{F_y}^{(\cT,\cG)}$ is weakly upper hemicontinuous in $L_p^\cG(T,X)$.

\item[F3] Suppose that $G$ is an $\cF$-measurable, $p$-integrably bounded and weakly compact valued correspondence from $T$ to $X$ with $1 \le p < \infty$. Let $F \colon T\times Y\to X$ be a weakly closed valued correspondence ($Y$ is a metric space). If
      \begin{enumerate}
        \item $F(t,y)\subseteq G(t)$ for $\lambda$-almost all $t\in T$;
        \item for any $y\in Y$, $F(\text{\ensuremath{\cdot}},y)$ is $\cF$-measurable;
        \item for any $t\in T$, $F(t,\cdot)$ is weakly lower hemicontinuous;
      \end{enumerate}
      then $H(y)= CI_{F_y}^{(\cT,\cG)}$ is weakly lower hemicontinuous in $L_p^\cG(T,X)$.
\end{enumerate}
\end{prop}

Before proving Proposition~\ref{prop-CI-hemicontinuity}, we first provide a simple lemma.

\begin{lem}\label{lem-convex upper}
Let $\{F_{\frac{1}{n}}\}_{n \ge 1} \cup \{F_0\}$ be a sequence of $\cT$-measurable correspondences from $T$ to $X$ such that $\{F_{\frac{1}{n}}(t)\}$ is weakly upper hemicontinuous at $0$ for $\lambda$-almost all $t\in T$. Suppose that $F_0$ is convex and weakly compact valued. Then the correspondence $\overline{\mbox{co}}\left(\cup_{m \ge n}\{F_{\frac{1}{m}}(t)\}\right)$ is also weakly upper hemicontinuous at $0$ for $\lambda$-almost all $t\in T$.
\end{lem}

\begin{proof}
It is obvious that the correspondence $\cup_{m \ge n}\{F_{\frac{1}{m}}(t)\}$ is weakly upper hemicontinuous at $0$ for $\lambda$-almost all $t\in T$. Since $F_0$ is convex and weakly compact valued, by Theorem~17.35 of \citet{AB2006}, the correspondence $\overline{\mbox{co}}\left(\cup_{m \ge n}\{F_{\frac{1}{m}}(t)\}\right)$ is also weakly upper hemicontinuous at $0$ for $\lambda$-almost all $t\in T$.
\end{proof}

Now we are ready to prove Proposition~\ref{prop-CI-hemicontinuity}.

\begin{proof}[Proof of Proposition~\ref{prop-CI-hemicontinuity}]
\

F1. The convexity is straightforward. We only need to show the weak compactness.

Similarly as that in the proof of Theorem~\ref{thm-ce-bochner}~(A2), it suffices to show that given an arbitrary sequence of $\cT$-measurable Bochner integrable selections $\{f_n\}_{n\in \bN}$ of $F$ with $g_n = E(f_n|\cG)$ for each $n\in \bN$, there is a subsequence of $\{g_n\}_{n\in \bN}$ that weakly converges in $L_p^\cG(T,X)$ to some point in $CI_{F}^{(\cT,\cG)}$. Following the same argument from the proof of Theorem~\ref{thm-ce-bochner}~(A2), the sequence $\{f_n\}_{n\in \bN}$ is relatively weakly compact in $L_p^\cT(T,X)$. Without loss of generality, we assume that $f_n$ weakly converges to some $f_0 \in L_p^\cT(T,X)$. Then $g_n$  weakly converges to $E(f_0|\cG)$ in $L_p^\cG(T,X)$.

By Theorem~29 in \citet[p.293]{RF2010}, there is a sequence $\{\phi_n\}$ that converges to $f_0$ in norm, and each $\phi_n$ is a convex combination of $\{f_n, f_{n+1}, \ldots\}$ for $n \ge 1$. Since $F$ is convex valued, $\phi_n$ is a selection of $F$, and hence $f_0$ is a selection of $F$.

\

F2. By Mazur's Theorem, $\overline{\mbox{co}}G$ is a convex and weakly compact valued correspondence. By (F1), $CI_{\overline{\mbox{co}}G}^{(\cT,\cG)}$ and $CI_{F_y}^{(\cT,\cG)}$ are weakly compact in $L_p^\cG(T,X)$. Note that when $\cG$ is a countably generated $\sigma$-algebra and $X$ is a separable Banach space, $L_p^\cG(T,X)$ is a separable Banach space. As shown in the proof of Theorem~\ref{thm-ce-bochner}~(A4), $CI_{F_y}^{(\cT,\cG)}$ can be viewed as a compact valued correspondence from $Y$ to a compact metric space. It suffices to show that if $y_{n}\rightarrow y_0$ in $Y$ (with $y_n \ne y_0$ for any $n \ge 1$), $g_n$ weakly converges to $g_0$ in $L_p^\cG(T,X)$ with $f_n$ being a $\cT$-measurable and Bochner integrable selection of $F_{y_n}$ and $g_n = E(f_n|\cG)$ for $n\ge 1$, then $g_0 \in CI_{F_{y_0}}^{(\cT,\cG)}$ (see Theorem~17.11 in \citet{AB2006}).

By the same argument as in (F1), we can assume without loss of generality that the sequence $\{f_n\}$ weakly converges to some $f_0 \in L_p^\cT(T,X)$, which implies that $g_0 = E(f_0|\cG)$. In addition, there is a sequence $\{\phi_n\}$ that converges to $f_0$ in norm, and each $\phi_n$ is a convex combination of $\{f_n, f_{n+1}, \ldots\}$ for $n \ge 1$. By Lemma~\ref{lem-convex upper}, $\overline{\mbox{co}}\left(\cup_{m \ge n}\{F_{y_m}(t)\}\right)$ is weakly upper hemicontinuous at $y_0$ for $\lambda$-almost all $t\in T$. Since $\phi_n$ is a $\cT$-measurable selection of $\overline{\mbox{co}}\left(\cup_{m \ge n}\{F_{y_m}(\cdot)\}\right)$, $f_0$ is a $\cT$-measurable selection of $F_{y_0}$.

\

F3. As in the proof of (F2) above, $CI_{F_y}^{(\cT,\cG)}$ can be viewed as a correspondence from $Y$ to the compact metric space $CI_{\overline{\mbox{co}}G}^{(\cT,\cG)}$. By Theorem~17.21 in \citet{AB2006}, it suffices to show that if $y_{n}\rightarrow y_0$ in $Y$ (with $y_n \ne y_0$ for any $n \ge 1$), then for any $g_0 \in CI_{F_{y_0}}^{(\cT,\cG)}$, there exists a subsequence of $\{y_{n_k}\}$ of $\{y_n\}$ and $g_{n_k} \in CI_{F_{y_{n_k}}}^{(\cT,\cG)}$ such that $g_{n_k}$ weakly converges to $g_0$ in $L_p^\cG(T,X)$.

Suppose that $f_0$ is a $\cT$-measurable selection of $F_{y_0}$ and $g_0 = E(f_0|\cG)$. Define a correspondence
$$\Phi(t) = \{(x_1, x_2, \ldots) \colon x_n \in F_{y_n} \mbox{ for each } n \ge 1, x_n \to f_0(t) \mbox{ weakly}\}.
$$
Since $F(t, \cdot)$ is weakly lower hemicontinuous, by Theorem~2 in \citet[p.27]{Hd1974}, $\Phi(t) \neq \emptyset$. Pick a $\cT$-measurable selection $\phi$ of $\Phi$. Then $\phi = (f_1, f_2, \ldots)$ with $f_n$ being a $\cT$-measurable selection of $F_{y_n}$ and $f_n(t)$ weakly converges to $f_0(t)$ for each $t$. Since $\{F(\cdot, y)\}$ is integrably bounded uniformly among all $y$, by dominated convergence theorem, $g_n = E(f_n|\cG)$ weakly converges to $E(f_0|\cG)$ in $L_p^\cG(T,X)$. This completes the proof.
\end{proof}

Proposition~\ref{prop-Bochener} is a simple corollary of Theorem~\ref{thm-ce-bochner}. Below, we prove Proposition~\ref{prop-norm integral}.

\begin{proof}[Proof of Proposition~\ref{prop-norm integral}]
\

B5. Pick a sequence of $\cT$-measurable Bochner integrable selections $\{f_n\}_{n\in \bN}$ of $F$. Let $x_n=\int_T f_n(t) \rmd\lambda(t)$ for each $n\in \bN$. We need to show that there is a subsequence of $\{x_n\}_{n\in \bN}$ which converges to some point in $I_{F}^{\cT}$.

We first show that $I_{F}^{\cT}$ is totally bounded and hence relatively compact.\footnote{The argument is the same as that in the proof of \citet[Theorem 2]{Sun1997}. We provide the details here for completeness.} Since $F$ is integrably bounded, there is a real valued integrable function $h$ on $(T,\cT,\lambda)$ such that $\sup\{\|x\|: x\in F(t)\} \leq h(t)$ for $\lambda$-almost all $t\in T$. Then for any $\epsilon > 0$, there is some $\delta > 0$ such that for any $\cT$-measurable set $D$ with $\lambda(D) < \delta$, $\int_D h \rmd \lambda < \frac{\epsilon}{2}$. By \citet[Proposition~3.2]{Sun1997}, there is a norm compact set $K_{\delta} \subseteq X$ such that
$$\lambda\{t \in T \colon F(t) \subseteq K_{\delta} \} > 1 - \delta.
$$
By Mazur's theorem, $\overline{\mbox{co}}K_{\delta}$ is a convex and norm compact set. Then one can choose finitely many points $z_1, \ldots, z_k \in \overline{\mbox{co}}K_{\delta}$ such that given any $x \in \overline{\mbox{co}}K_{\delta}$, $\|x - z_i\| < \frac{\epsilon}{2}$ for some $1 \le i \le k$. Pick any $\cT$-measurable selection $f$ of $F$. Let $D_{\delta} = f^{-1}(\overline{\mbox{co}}K_{\delta})$. Then $\lambda(T\setminus D_{\delta}) < \delta$, which implies that
$$\|\int_{T\setminus D_{\delta}} f \rmd \lambda\| < \int_{T\setminus D_{\delta}} h \rmd \lambda < \frac{\epsilon}{2}.
$$
Since $\int_{D_{\delta}} f \rmd \lambda \in \overline{\mbox{co}}K_{\delta}$, $\|\int_{D_{\delta}} f \rmd \lambda - z_i\| < \frac{\epsilon}{2}$ for some $1 \le i \le k$, and hence $\|\int_{T} f \rmd \lambda - z_i\| < \epsilon$. Thus, $I_{F}^{\cT}$ is totally bounded.

Without loss of generality, we assume that $\{x_n\}_{n\in \bN}$ converges to some point $x_0\in X$. By Lemma~\ref{lem-distribution} (2), $D_F^\cT$ is compact. Let $\mu_n=\lambda \circ f_n^{-1}$. There is a subsequence of $\{\mu_n\}_{n\in \bN}$, say itself, which weakly converges to some $\mu_0\in  D_F^\cT \subseteq \cM(X)$. Suppose that $\mu_0=\lambda \circ f_0^{-1}$, where $f_0$ is a Bochner integrable selection of $F$. Given  any continuous linear functional $x^*$ on $X$, since $x^*\circ f_n$ is uniformly integrable and $\lambda \circ f_n^{-1}\circ (x^*)^{-1}$ weakly converges to $\lambda \circ f_0^{-1}\circ (x^*)^{-1}$, by Theorem 5.4 of \citet[p.32]{Billingsley1968}, we have that
$$\int_T x^*\circ f_n \rmd\lambda(t)\to \int_T x^*\circ f_0 \rmd\lambda(t).$$
Thus, $x^*(x_0)=x^*(\int_T f_0(t) \rmd\lambda(t))$, which implies that $x_0=\int_T f_0(t) \rmd\lambda(t)$. Therefore, $I_{F}^{\cT}$ is norm compact.

\

B6. By (B5), we know that $I_{F_{y}}^{\cT}$ is a norm compact valued correspondence from $Y$ to the compact space $I_{G}^{\cT}$. Let $\{x_n\}_{n=0}^\infty$ be a sequence in $X$ and $\{y_n\}_{n=0}^\infty$ be a sequence in $Y$ ($y_n \neq y_0$ for any $n \ge 1$). Suppose that $x_n \in I_{F_{y_n}}^{\cT}$ for each $n\ge 1$, $\lim_{n\to\infty} x_n = x_0$ and $\lim_{n\to\infty} y_n = y_0$. We need to show that $x_0 \in I_{F_{y_0}}^{\cT}$.

It is clear that $x_0 \in \limsup_{n\to\infty} I_{F_{y_n}}^{\cT}$.  For each $t\in T$, since the correspondence $F(t, \cdot)$ is norm upper hemicontinuous,  $\limsup_{n\to\infty} F(t, y_n) \subseteq F(t, y_0)$. Let $\{A_n\}_{n\in\bN}$ be a sequence of sets in a compact metric space, and $A_0$ be a compact set. By Theorem 1 in \citet[p.17]{Hd1974}, $\limsup_{n\to\infty} A_n \subseteq A_0$ if and only if the Hausdorff semi-distance $\sigma(A_n,A_0)$ goes to $0$ as $n\to \infty$. Since $F_{y_0}(t)$ is norm compact, we have that for each $t \in T$,
$$\lim_{n\to\infty} \sigma (F_{y_n}(t), F_{y_0}(t)) =0.
$$
By Proposition 4.3 of \citet{Sun1997},  $\lim_{n\to\infty} \sigma (I_{F_{y_n}}^{\cT}, I_{F_{y_0}}^{\cT}) = 0$. By (B5), $I_{F_{y_0}}^{\cT}$ is norm compact, and hence $\limsup_{n\to\infty} I_{F_{y_n}}^{\cT} \subseteq I_{F_{y_0}}^{\cT}$. This implies that $x_0 \in \limsup_{n\to\infty} I_{F_{y_n}}^{\cT} \subseteq I_{F_{y_0}}^{\cT}$.
\end{proof}

The following result follows from Proposition~\ref{prop-norm integral}~(B6).

\begin{coro}\label{coro-norm upper}
Suppose that $G$ is an $\cF$-measurable, integrably bounded and norm compact valued correspondence from $T$ to $X$. Let $\{F_n\}^\infty_{n=1}$ be a sequence of $\cF$-measurable correspondences from $T$ to $X$ such that $F_n(t)\subseteq G(t)$ for $\lambda$-almost all $t\in T$ and $n\geq 1$. Let $F_0 =\limsup_{n \to \infty} F_n$. Then $\limsup_{n\to\infty} I_{F_n}^{\cT} \subseteq I_{F_0}^{\cT}$.
\end{coro}

\subsection{Proofs of the results in Section~\ref{sec-G integral}}\label{proof-G}

The following lemma is analogous to Lemma~\ref{lem-separate}. We skip its proof as it follows the same argument as that in the proof of Lemma~\ref{lem-separate}.

\begin{lem}\label{lem-separate_2}
Suppose that $\tilde{\Phi}$ is the set of all mappings $\varphi$ from $(T,\cG,\lambda)$ to $X$ such that $\varphi(t)\in \{x_m\}_{m \ge 1}$ for $\lambda$-almost all $t\in T$. Then the set $\tilde{\Phi}$ separates points in $L_p^\cG(T,X^*)$ for $1 < p \le \infty$.
\end{lem}

Recall that $X$ is separable. The following result has been shown in \citet[p.149]{Sun1997}. We summarize it as Lemma~\ref{lem-sun1997-2} and provide the proof for completion.

\begin{lem}\label{lem-sun1997-2}
There is a countable dense set $\{x_m\}_{m \ge 1}$ in the unit ball of $X$ such that
\begin{enumerate}
\item the weak$^*$ topology coincides with the topology induced by the metric $d_w$ on any norm bounded subset of $X^*$, where for any $x^*$ and $y^*$ in $X^*$,
$$d_w(x^*,y^*) = \sum^\infty_{m=1} {1\over{2^m}} |x^*(x_m) - y^*(x_m)|;
$$

\item the Borel $\sigma$-algebras on $X^*$ generated by the weak$^*$~topology and the topology induced by the metric $d_w$ are the same.
\end{enumerate}
Denote the completion of $(X^*, d_w)$ by $(Z^*,d_w)$.
\end{lem}

\begin{proof}
(1) Since X is separable, there exists a sequence $\{x_m\}_{m \ge 1}$ in the unit ball of $X$ such that the linear space spanned by this sequence is dense in $X$. Define the metric $d_w$ on $X^*$ as described in the lemma. The topology induced by the metric $d_w$ is weaker than the weak$^*$ topology on $X^*$, which means that on any weak$^*$ compact subset, and hence on any norm bounded subset of $X^*$, these two topologies coincide.

(2) As $X^*$ is a union of countably many weak$^*$ compact subsets, it is a Borel set in $(Z^*,d_w)$. Then the Borel $\sigma$-algebras on $X^*$ generated by the weak$^*$ topology and the topology induced by the metric $d_w$ are the same.
\end{proof}

Given a Gel$^\prime$fand integrable mapping $f$ from $(T, \cT, \lambda)$ to $X^*$ (endowed with the weak$^*$ topology), by Lemma~\ref{lem-sun1997-2} above, one can also view $f$ as a measurable mapping from $(T, \cT, \lambda)$ to $(Z^*, d_w)$. The following lemma is analogous to Lemma~\ref{lem-Bochner ce}. We skip its proof as it is the same as that in the proof of Lemma~\ref{lem-Bochner ce}.

\begin{lem}\label{lem-G ce}
For a Gel$^\prime$fand integrable mapping $f$ from $(T, \cT, \lambda)$ to $X^*$,
$$E(f|\cG)(t) = \int_{X^*} I_{X^*}(x^*) \lambda^{f|\cG}(t,\rmd x^*)
$$
for $\lambda$-almost all $t\in T$, where (1) $I_{X^*}(\cdot)$ is the identity mapping from $X^*$ (endowed with the topology induced by the metric $d_w$) to $X^*$ (endowed with the weak$^*$ topology), and (2) $\lambda^{f|\cG}$ is the regular conditional distribution given $\cG$ generated by viewing $f$ as a mapping from $(T, \cT, \lambda)$ to $(Z^*, d_w)$.
\end{lem}

We first prove the following result, which parallels Proposition~\ref{prop-CI-hemicontinuity} in terms of Gel$^\prime$fand integral and conditional expectation. This proposition will be useful for proving Theorem~\ref{thm-ce-G}. Compared with Theorem~\ref{thm-ce-G}, we impose the convexity condition on the correspondence $F$ in Proposition~\ref{prop-CG-hemicontinuity}~(G1) and (G2), and show that the results on compactness and preservation of upper hemicontinuity still hold without the nowhere equivalence condition. In (G3), we prove the preservation of lower hemicontinuity without either the nowhere equivalence condition, or the convexity condition.

\begin{prop}\label{prop-CG-hemicontinuity}
For any sub-$\sigma$-algebra $\cG$ of $\cF$, we have the following properties.
\begin{enumerate}
\item[G1] For any $\cF$-measurable, $p$-integrably bounded, convex and weak$^*$~compact valued correspondence $F$ from $T$ to $X^*$ with $1 < p \leq \infty$, $CG_{F}^{(\cT,\cG)}$ is convex and weak$^*$~compact (in the dual space of $L_q^\cG(T,X)$ with $\frac{1}{p} + \frac{1}{q} = 1$).

\item[G2] Suppose that $G$ is an $\cF$-measurable, $p$-integrably bounded, and weak$^*$ compact valued correspondence from $T$ to $X^*$ with $1< p \le \infty$. Let $F$ be a convex and weak$^*$ closed valued correspondence from $T\times Y\to X^*$ ($Y$ is a metric space). If
      \begin{enumerate}
        \item $F(t,y) \subseteq G(t)$ for $\lambda$-almost all $t\in T$;
        \item for any $y\in Y$, $F(\text{\ensuremath{\cdot}},y)$ is $\cF$-measurable;
        \item for any $t\in T$, $F(t,\cdot)$ is weak$^*$ upper hemicontinuous;
      \end{enumerate}
      then $H(y)= CG_{F_y}^{(\cT,\cG)}$ is weak$^*$ upper hemicontinuous (in the dual space of $L_q^\cG(T,X)$).

\item[G3] Suppose that $G$ is an $\cF$-measurable, $p$-integrably bounded, and weak$^*$ compact valued correspondence from $T$ to $X^*$ with $1< p \le \infty$. Let $F$ be a weak$^*$ closed valued correspondence from $T\times Y$ to $X^*$ ($Y$ is a metric space). If
      \begin{enumerate}
        \item $F(t,y) \subseteq G(t)$ for $\lambda$-almost all $t\in T$;
        \item for any $y\in Y$, $F(\text{\ensuremath{\cdot}},y)$ is $\cF$-measurable;
        \item for any $t\in T$, $F(t,\cdot)$ is weak$^*$ lower hemicontinuous;
      \end{enumerate}
      then $H(y)= CG_{F_y}^{(\cT,\cG)}$ is weak$^*$ lower hemicontinuous (in the dual space of $L_q^\cG(T,X)$).
\end{enumerate}
\end{prop}

\begin{proof}

Recall that $\cF$ is countably generated. Fix a sub-$\sigma$-algebra $\cG$ of $\cF$.

\

G1. It is obvious that $CG_{F}^{(\cT,\cG)}$ is convex. We only need to show the compactness.

Since $F$ is $p$-integrably bounded, $CG_{F}^{(\cT,\cG)}$ is norm bounded in $L_p^\cG(T,X^*)$. By Alaoglu's Theorem (see Theorem~11.7.1 in \citet{Loeb2016}), the closed unit ball of the dual space of $L_q^\cG(T,X)$ is weak$^*$ compact with $\frac{1}{p} + \frac{1}{q} = 1$.\footnote{Note that $L_p^\cG(T,X^*)$ is only a subset of the dual space of $L_q^\cG(T,X)$, as we do not assume the Radon-Nikod\'ym property on the space $X^*$; see Theorem 1 in \citet[p.98]{DU1977}.} Since $X$ is separable and $\cG$ is countably generated, $L_q^\cG(T,X)$ is a separable Banach space. Then the weak$^*$ topology on the closed unit ball of the dual space of $L_q^\cG(T,X)$ is metrizable (see Corollary~11 in \citet[p.306]{RF2010}). Thus, we can work with sequence instead of net. Let $\{f_n\}_{n \ge 1}$ be a sequence of $\cT$-measurable Gel$^\prime$fand integrable selections of $F$. Denote $\cF_1 = \sigma(h, \{f_n\}, \cF\}$, the smallest $\sigma$-algebra such that $\cF \subseteq \cF_1$, and $h$ and $\{f_n\}$ are $\cF_1$-measurable. Then $\cF_1$ is a countably generated $\sigma$-algebra.

%Following the same argument as before, the closed unit ball of the dual space of $L_q^{\cF_1}(T,X)$ is weak$^*$ compact. Thus, $\{f_n\}$ has a weak$^*$ convergent subsequence. Without loss of generality, we assume that $\{f_n\}$ weak$^*$ converges to some $f_{0}$ in the dual space of $L_q^{\cF_1}(T,X)$. We shall show that $f_{0} \in CG_{F}^{(\cT,\cF_1)}$.

Let $\eta$ be the uniform distribution on the unit interval $I = [0,1]$, endowed with the Borel $\sigma$-algebra $\cB$. Denote $T' = T \times I$, $\cT' = \cT \otimes \cB$, $\lambda' = \lambda \otimes \eta$, and $\cF'_1 = \cF_1 \otimes \{I, \emptyset\}$. Define $F'$ as an $\cF'_1$-measurable correspondence such that $F'(t, i) = F(t)$ for any $(t, i) \in T'$. Similarly, we define a sequence $\{f'_n\}$ such that $f'_n(t, i) = f_n(t)$ for $(t, i) \in T'$ and $n \ge 1$. Then $f'_n$ is an $\cF'_1$-measurable selection of $F'$. By Lemma~\ref{lem-sun1997-2}, one can also view $F'$ as a compact valued correspondence from $(T', \cF'_1, \lambda')$ to $(Z^*, d_w)$, taking values from $X^*$. Then $\cT'$ is nowhere equivalent to $\cF'_1$ under $\lambda'$. By Lemma~\ref{lem-distribution}, $\cR^{(\cT', \cF'_1)}_{F'}$ is weakly compact, which implies that there exists some $\cT'$-measurable selection $f'$ of $F'$ such that a subsequence of $\lambda'^{f'_n|\cF'_1}$, say itself, converges to $\lambda'^{f'|\cF'_1}$ weakly. Again by Lemma~\ref{lem-sun1997-2}, $f'$ can be viewed as a $\cT'$-measurable mapping from $T'$ to $X^*$ endowed with the weak$^*$ topology.

Let $f'_{\infty} = E(f' | \cF'_1)$. Then $f'_{\infty}$ is an $\cF'_1$-measurable mapping from $T'$ to $X^*$ (endowed with the weak$^*$ topology). Since $F$ is convex valued, so is $F'$. Thus, $f'_{\infty}$ is also a selection of $F'$. We shall show that $\{f'_n\}$ weak$^*$ converges to $f'_{\infty}$ in $L_p^{\cF'_1}(T', X^*)$.

Recall that in Lemma~\ref{lem-sun1997-2}, the metric $d_w$ is constructed by the countable dense set $\{x_m\}_{m \ge 1}$. For each $x_m$ and a nontrivial subset $E \in \cF'_1$ with $\lambda'(E) > 0$,
\begin{align*}
\int_{T'} 1_E(t') f'_n(t', x_m) \rmd \lambda'(t')
& = \int_{T'} 1_E(t') E(f'_n | \cF'_1)(t', x_m) \rmd \lambda'(t') \\
& = \int_{T'} \int_{X^*} 1_E(t') I_{X^*}(x^*, x_m) \lambda'^{f'_n|\cF'_1} (t',\rmd x^*) \rmd \lambda'(t') \\
& \to \int_{T'} \int_{X^*} 1_E(t') I_{X^*}(x^*, x_m) \lambda'^{f'|\cF'_1} (t',\rmd x^*) \rmd \lambda'(t') \\
& = \int_{T'} 1_E(t') E(f' | \cF'_1)(t', x_m) \rmd \lambda'(t') \\
& = \int_{T'} 1_E(t') f'_{\infty} (t', x_m) \rmd \lambda'(t').
\end{align*}
The first equality holds since $f'_n$ is $\cF'_1$-measurable. The second and third equalities are due to Lemma~\ref{lem-G ce}. The convergence is true since $\lambda'^{f'_n|\cF'_1}$ converges to $\lambda'^{f'|\cF'_1}$ weakly. The last equality is by the definition of $f'_{\infty}$.

Let $\Phi_1$ be the set of all mappings $\varphi_1$ from $(T', \cF'_1, \lambda')$ to $X$ such that $\varphi_1(t') \in \{x_m\}_{m \ge 1}$ for $\lambda'$-almost all $t' \in T'$. Given $\varphi_1 \in \Phi_1$, let $E_m = \{t' \in T' \colon \varphi_1(t') = x_m \}$ for $m \ge 1$. Then as $n \to \infty$,
\begin{align*}
\int f'_n \varphi_1 \rmd \lambda'
& = \sum_{m \ge 1} \int_{T'} 1_{E_m}(t') f'_n(t', x_m) \rmd \lambda'(t') \\
& \to \sum_{m \ge 1} \int_{T'} 1_{E_m}(t') f'_{\infty}(t', x_m) \rmd \lambda'(t') \\
& = \int f'_{\infty} \varphi_1 \rmd \lambda'.
\end{align*}

Suppose that $\{f'_n\}$ does not weak$^*$ converge to $f'_{\infty}$ in $L_p^{\cF'_1}(T', X^*)$. Then there exists a subsequence of $\{f'_n\}$, say itself, some $\psi \in L_q^{\cF'_1}(T', X)$, and $\tilde{\epsilon} > 0$ such that
$$\int f'_n \psi \rmd \lambda'  - \int f'_{\infty} \psi \rmd \lambda' > \tilde{\epsilon}.
$$
Without loss of generality, we assume that $\psi \in L_{\infty}^{\cF'_1}(T', X)$ and is essentially bounded by $1$.

For $m \ge 1$, let
$$J_m = \left\{ t' \in T' \colon \|x_m - \psi(t')\|_X \le \min\{1, \frac{\tilde{\epsilon}}{3 h_t}\} \right\},
$$
where $\|\cdot\|_X$ is the norm on the Banach space $X$. Define $\psi_1(t') = x_m$ for $t' \in J_m \setminus \cup_{0 \le k < m} J_k$, where $J_0 = \emptyset$. For each $n$,
\begin{align*}
\int_{T'} f'_n(t') \psi_1(t') \rmd \lambda'(t')  - \int_{T'} f'_n(t') \psi(t') \rmd \lambda'(t')
& = \int_{T'} f'_n(t') (\psi_1(t') - \psi(t') ) \rmd \lambda'(t') \\
& \ge - \int_{T'} |f'_n(t') (\psi_1(t') - \psi(t') )| \rmd \lambda'(t') \\
& = - \int_{T'} \|f'_n(t')\| \|(\psi_1(t') - \psi(t') )\|_X \rmd \lambda'(t') \\
& \ge - \int_{T'} h(t') \cdot \min\{1, \frac{\tilde{\epsilon}}{3 h(t')}\} \rmd \lambda'(t') \\
& \ge - \frac{\tilde{\epsilon}}{3}.
\end{align*}
Similarly,
$$\int_{T'} f'_{\infty}(t') \psi(t') \rmd \lambda'(t')  - \int_{T'} f'_{\infty}(t') \psi_1(t') \rmd \lambda'(t') \ge - \frac{\tilde{\epsilon}}{3}.
$$
Then for each $n$,
\begin{align*}
\int f'_n \psi_1 \rmd \lambda'  - \int f'_{\infty} \psi_1 \rmd \lambda'
& = \left( \int f'_n \psi_1 \rmd \lambda'  - \int f'_{n} \psi \rmd \lambda' \right) \\
& + \left( \int f'_n \psi \rmd \lambda'  - \int f'_{\infty} \psi \rmd \lambda' \right) \\
& + \left( \int f'_{\infty} \psi \rmd \lambda'  - \int f'_{\infty} \psi_1 \rmd \lambda' \right) \\
& \ge \frac{\tilde{\epsilon}}{3}.
\end{align*}
This is a contradiction. Thus, $\{f'_n\}$ weak$^*$ converges to $f'_{\infty}$ in $L_p^{\cF'_1}(T', X^*)$.

Since $f'_{\infty}$ is an $\cF'_1$-measurable selection of $F'$, there exists some $\cF_1$-measurable selection $f_{\infty}$ of $F$ such that $f'_{\infty}(t, i) = f_{\infty}(t)$ for $\lambda'$-almost all $(t, i)$. It is clear that $f_n$ weak$^*$ converges to $f_{\infty}$ in $L_p^{\cF_1}(T, X^*)$. Let $g_n = E(f_n|\cG)$ for each $n \ge 1$. Then $g_n$ weak$^*$ converges to $g_{\infty} = E(f_{\infty}|\cG)$ in $L_p^{\cG}(T,X^*)$. This completes the proof.

\

G2. Since $G$ is $p$-integrably bounded, $CG_{w^*-\overline{\mbox{co}}G}^{(\cT,\cG)}$ is norm bounded in $L_p^\cG(T,X^*)$. Note that as $\cG$ is countably generated and $X$ is separable, $L_q^\cG(T,X)$ is a separable Banach space with $\frac{1}{p} + \frac{1}{q} = 1$. By \citet[Corollary~11, p.306]{RF2010}, the weak$^*$ topology on the norm bounded closed ball of the dual space of $L_q^\cG(T,X)$ is metrizable. By (G1), $CG_{w^*-\overline{\mbox{co}}G}^{(\cT,\cG)}$ and $CG_{F_y}^{(\cT,\cG)}$ are weak$^*$ compact in the dual space of $L_q^\cG(T,X)$. Then $CG_{F_y}^{(\cT,\cG)}$ can be viewed as a compact valued correspondence from $Y$ to a compact metric space. It suffices to show that if $y_{n}\rightarrow y_{\infty}$ in $Y$ (with $y_n \ne y_{\infty}$ for $n \ge 1$), $g_n$ weak$^*$ converges to $g_{\infty}$ in $L_p^\cG(T,X^*)$ with $f_n$ being a $\cT$-measurable and Gel$^\prime$fand integrable selection of $F_{y_n}$ and $g_n = E(f_n|\cG)$ for $n \ge 1$, then $g_{\infty} \in CG_{F_{y_{\infty}}}^{(\cT,\cG)}$ (see Theorem~17.11 in \citet{AB2006}).

By (G1), there is a subsequence of $\{f_n\}$, say itself, weak$^*$ converging to some $f_{\infty} \in L_p^{\cF_1}(T,X^*)$, where $\cF_1$ is defined in (G1) and $f_{\infty}(t) \in w^*-\overline{\mbox{co}}(\cup_{k \ge n}\{f_k(t)\})$ for any $n \ge 1$. It implies that $g_{\infty} = E(f_{\infty}|\cG)$. Fix any $t\in T$ such that $F_{y_n}(t)$ is weak$^*$ upper hemicontinuous at $y_{\infty}$ as $y_n \to y_{\infty}$. By \citet[Theorem~17.35]{AB2006}, $w^*-\overline{\mbox{co}}\left(\cup_{k \ge n}\{F_{y_k}(t)\}\right)$ is also weak$^*$ upper hemicontinuous as $n \to \infty$. Thus, $f_{\infty}$ is a selection of $F_{y_{\infty}}$. This completes the proof.

\

G3. The proof is the same as that for Proposition~\ref{prop-CI-hemicontinuity}~(F3).
%
%We provide the proof for completeness.
%
%
%As in the proof of (G2) above, $CG_{F_y}^{(\cT,\cG)}$ can be viewed as a correspondence from $Y$ to the compact metric space $CG_{w^*-\overline{\mbox{co}}G}^{(\cT,\cG)}$. By Theorem~17.21 in \citet{AB2006}, it suffices to show that if $y_{n}\rightarrow y_0$ in $Y$ (with $y_n \ne y_0$ for any $n \ge 1$), then for any $g_0 \in CG_{F_{y_0}}^{(\cT,\cG)}$, there exists a subsequence of $\{y_{n_k}\}$ of $\{y_n\}$ and $g_{n_k} \in CG_{F_{y_{n_k}}}^{(\cT,\cG)}$ such that $g_{n_k}$ weak$^*$ converges to $g_0$ in $L_p^\cG(T,X^*)$.
%
%
%Suppose that $f_0$ is a $\cT$-measurable selection of $F_{y_0}$ and $g_0 = E(f_0|\cG)$. Define a correspondence
%$$\Phi(t) = \{(x_1, x_2, \ldots) \colon x_n \in F_{y_n} \mbox{ for each } n \ge 1, x_n \mbox{ weak}^* \mbox{ converges to } f_0(t) \}.
%$$
%Since $F(t, \cdot)$ is weak$^*$ lower hemicontinuous, by Theorem~2 in \citet[p.27]{Hd1974}, $\Phi(t) \neq \emptyset$. Pick a $\cT$-measurable selection $\phi$ of $\Phi$. Then $\phi = (f_1, f_2, \ldots)$ with $f_n$ being a $\cT$-measurable selection of $F_{y_n}$, and $f_n(t)$ weak$^*$ converges to $f_0(t)$ for each $t$. Since $\{F(\cdot, y)\}$ is integrably bounded uniformly among all $y$, by dominated convergence theorem, $g_n = E(f_n|\cG)$ weak$^*$ converges to $E(f_0|\cG)$ in $L_p^\cG(T,X^*)$. This completes the proof.
\end{proof}

Following the same argument as in Remark~\ref{rmk-distribution}, one can obtain the desirable properties for regular conditional distribution of correspondences based on Proposition~\ref{prop-CG-hemicontinuity}. We summarize the results in the corollary below.

\begin{coro}\label{coro-convex distribution}
Let $\cF$ be a countably-generated sub-$\sigma$-algebra of $\cT$, and $X$ a Polish space. For any $\sigma$-algebra $\cG \subseteq \cF$, we have the following properties.
\begin{enumerate}
%  \item For any convex and closed valued $\cF$-measurable correspondence $F$ from $T$ to $X$, $\cR^{(\cT,\cG)}_F$ is convex.

  \item For any convex and compact valued $\cF$-measurable correspondence $F$ from $T$ to $X$, $\cR^{(\cT,\cG)}_F$ is convex and weakly compact.

  \item Let $F$ be a compact valued $\cF$-measurable correspondence from $T$ to $X$, and $G$ a convex and closed valued correspondence from $T\times Y$ to $X$ ($Y$ is a metric space) such that
      \begin{enumerate}
        \item[a] for any $(t,y)\in T\times Y$, $G(t,y)\subseteq F(t)$;
        \item[b] for any $y\in Y$, $G(\text{\ensuremath{\cdot}},y)$ (denoted as $G_y$) is $\cF$-measurable from $T$ to $X$;
        \item[c] for any $t\in T$, $G(t,\cdot)$ (denoted as $G_t$) is upper hemicontinuous from $Y$ to $X$.
      \end{enumerate}
      The correspondence $H(y)=\cR_{G_{y}}^{(\cT,\cG)}$ is upper hemicontinuous from $Y$ to $\cR^{\cG}$.

      \item Let $F$ be a compact valued $\cF$-measurable correspondence from $T$ to $X$, and $G$ a closed valued correspondence from $T\times Y$ to $X$ ($Y$ is a metric space) such that
      \begin{enumerate}
        \item[a] for any $(t,y)\in T\times Y$, $G(t,y)\subseteq F(t)$;
        \item[b] for any $y\in Y$, $G(\text{\ensuremath{\cdot}},y)$ (denoted as $G_y$) is $\cF$-measurable from $T$ to $X$;
        \item[c] for any $t\in T$, $G(t,\cdot)$ (denoted as $G_t$) is lower hemicontinuous from $Y$ to $X$.
      \end{enumerate}
      The correspondence $H(y)=\cR_{G_{y}}^{(\cT,\cG)}$ is lower hemicontinuous from $Y$ to $\cR^{\cG}$.
  \end{enumerate}
\end{coro}

Now we are ready to prove Theorem~\ref{thm-ce-G}.

\begin{proof}[Proof of Theorem~\ref{thm-ce-G}]

%Given a sub-$\sigma$-algebra $\cG$, suppose that $\psi$  is a mapping from $T$ to some Polish space $Y$ such that $\cG$ is induced by $\psi$,  $\eta$=$\lambda %\psi^{-1}$.
%Then $F_1=\{\psi,F\}$ is an $\cF$-measurable correspondence from $T$ to $Y\times X$.
%Let $g=(\psi,f)$ be a $\cT$-measurable selection of $F_1$ and $\nu^{(\psi,f)}=\lambda g^{-1}$. Then $\nu^{(\psi,f)}_Y=\eta$.
%Since $X$ and $Y$ are both Polish spaces, there exists a Borel family of probability measures $\{\nu^{(\psi,f)}(y,\cdot)\}_{y\in Y}$ (which is $\eta$-a.e. %uniquely determined) in $\cM(X)$ such that
%$$\nu^{(\psi,f)}=\int_Y \nu^{(\psi,f)}(y,\cdot)\rmd\eta(y).
%$$
%Let $f^\psi(t) = \int_X  \rmd \nu^{(\psi,f)}(\psi(t), x)$ for each $t\in T$. Then $f^\psi$ is $\cG$-measurable mapping from $T$ to $X$. Let
%$$\cE_\psi=\{f^\psi\colon f \mbox{ is a } \cT\mbox{-measurable Bochner integrable selection of } F\},
%$$
%we shall prove that $\cE_\psi$ coincides with $E_F^{(\cT,\cG)}$.

\

C1. The proof is the same as that in the proof of Theorem~\ref{thm-ce-bochner}~(A1).

%The correspondence $F$ is $\cF$-measurable from $T$ to the Polish space $(Z^*, d_w)$. Given two $\cT$-measurable Bochner integrable selections $f_1,f_2$ of $F$ and $0\leq \alpha \leq 1$. Since $\cR^{(\cT,\cG)}_F$ is convex, there exists a $\cT$-measurable selection $f$ of $F$ such that $\mu^{f|\cG} = \alpha \mu^{f_1|\cG} + (1-\alpha) \mu^{f_2|\cG}$. Then for $\lambda$-almost all $t\in T$,
%\begin{align*}
%E(f|\cG)(t)
%& = \int_X I_X \mu^{f|\cG}(t,\rmd x) = \int_X I_X \left(\alpha \mu^{f_1|\cG}(t, \rmd x) + (1-\alpha) \mu^{f_2|\cG}(t, \rmd x) \right)\\
%& = \alpha E(f_1|\cG)(t) + (1-\alpha) E(f_2|\cG)(t).
%\end{align*}
%That is, $CI_{F}^{(\cT,\cG)}$ is convex.

\

C2. By the proof of Proposition~\ref{prop-CG-hemicontinuity}~(G1), $CG_{w^*-\overline{\mbox{co}}F}^{(\cT,\cG)}$ is metrizable and weak$^*$ compact. It suffices to show that given any sequence of $\cT$-measurable and Gel$^\prime$fand integrable selections $\{f_n\}$ of $F$ with $g_n = E(f_n|\cG)$ for $n \ge 1$, there is a subsequence of $\{g_n\}$, say itself, which weak$^*$ converges to some point in $CG_{F}^{(\cT,\cG)}$.

By Proposition~\ref{prop-CG-hemicontinuity}~(G1) again, a subsequence of $\{f_n\}$, say itself, weak$^*$ converges to some $f_0$ in $L_p^{\cG_1}(T, X^*) \subseteq L_p^{\cT}(T, X^*)$. Thus, $g_n$ weak$^*$ converges to $E(f_0|\cG)$ in $L_p^{\cG}(T, X^*)$. Denote $g_0 = E(f_0|\cG)$.

Since the correspondence $F$ is compact valued from $(T, \cF,\lambda)$ to the Polish space $(Z^*, d_w)$,  $\cR^{(\cT,\cG)}_F$ is compact by Lemma~\ref{lem-distribution}. Thus, there is a $\cT$-measurable Gel$^\prime$fand integrable selection $f$ of $F$ and a subsequence of $\{f_{n}\}_{n\in \bN}$, say itself,  such that $\mu^{f_n|\cG}$ weakly converges to $\mu^{f|\cG}$ as $n \to \infty$. Given any $\varphi$ from $(T,\cG,\lambda)$ to $X$ such that $\varphi(t)\in \{x_m\}_{m\in \bN}$ for $\lambda$-almost all $t\in T$. Repeating the argument in the proof of Theorem~\ref{thm-ce-bochner}~(A2), we could show that
$$\lim_{n\to\infty} \int_T \varphi^{**}(t,g_n(t)) \rmd\lambda  = \int_T \varphi^{**}(t,g(t)) \rmd\lambda,$$
where $g = E(f|\cG)$. By Lemma~\ref{lem-separate_2}, $g_0 = g \in CG_{F}^{(\cT,\cG)}$, and hence $CG_{F}^{(\cT,\cG)}$ is weak$^*$ compact.

\

C3. We prove by contradiction. Suppose that there is a $\cT$-measurable Gel$^\prime$fand integrable selection $f$ of $w^* - {\overline {co}}F$ such that $g= E(f|\cG) \notin CG^{(\cT,\cG)}_F$.

By (C1) and (C2), $CG_{F}^{(\cT,\cG)}$ is a convex and weak$^*$ compact set in the dual space of $L_q^\cG(T,X)$. Moreover, $L^\cG_q(T,X)$ is separable with $\frac{1}{p}+\frac{1}{q} = 1$. Then the weak$^*$ topology on the norm bounded subset of $L^\cG_p(T,X^*)$ is metrizable (see the proof of (G1)). We choose a countable dense set $\{\varphi_m\}_{m\in\bN}$ in the unit ball of $L^\cG_q(T,X)$ and define the metric $u_w$, where
$$u_w(g_1^*,g_2^*) = \sum^\infty_{m=1} {1\over{2^m}} |g_1^*(\varphi_m) - g_2^*(\varphi_m)|
$$
for each pair of $g_1^*$ and $g_2^*$ in $L^\cG_p(T,X^*)$. Let $(L^*,u_w)$ be the completion of $(L^\cG_p(T,X^*), u_w)$.

Since $g \notin CG^{(\cT,\cG)}_F$ and $CG_{F}^{(\cT,\cG)}$ is convex and compact in $(L^*,u_w)$,
$$u_w(g, CG_{F}^{(\cT,\cG)}) = \min \{u_w(g, g') \colon g' \in CG_{F}^{(\cT,\cG)}\}> 0.
$$
We choose $M > 0$ and sufficiently large $n$ such that
\begin{enumerate}
\item $CG_{F}^{(\cT,\cG)}$ is norm bounded by $M$,

\item $u_w(g, CG_{F}^{(\cT,\cG)}) > \frac{3M}{2^n}$.
\end{enumerate}
Then for any $g'\in CG_{F}^{(\cT,\cG)}$,
$$\sum_{m=1}^n \frac{1}{2^m} |\varphi^{**}_m (g - g')| > \frac{3M}{2^n} - \sum_{m=n+1}^\infty \frac{2M}{2^m} = \frac{M}{2^n}.
$$

Define a linear mapping $\gamma$ from $L^\cG_p(T,X^*)$ to $\bR^n$ as $\gamma(g) = (\varphi_1^{**}(g), \cdots, \varphi_n^{**}(g))$ for any $g\in L^\cG_p(T,X^*)$. Then $\gamma(g)$ is not in the compact set $\gamma\big(CG_{F}^{(\cT,\cG)}\big)$. Define a new correspondence $F_1$ from $T$ to $\bR^n$:
$$F_1(t) = \left\{ (a_1,\ldots,a_n) \colon a_i = x^*(\varphi_i(t)), i=1,\ldots,n,  x^* \in F(t) \right\}.
$$
Then $F_1 = \gamma \circ F$ is $\cF$-measurable. By Filippov's implicit function theorem,  any measurable selection of $F_1$ is the composition of $\gamma$ with a measurable selection of $F$. Then $\gamma\big(CG_{F}^{(\cT,\cG)}\big) = G_{F_1}^{\cT}$. Since $\gamma$ is  continuous, $F_1$ is compact valued in $\bR^n$. It is clear that  $\gamma(w^* -{\overline {co}}F(t)) = {\overline {co}} F_1(t) = co (F_1(t))$. We have that
$$\gamma(g) =   \gamma(f)  \in G^\cT_{co(F_1)} = G^\cT_{F_1} = \gamma\big(CG_{F}^{(\cT,\cG)}\big),
$$
where the equality $G^\cT_{co(F_1)} = G^\cT_{F_1}$ is due to Theorem~4 of \citet[p.64]{Hd1974}. We arrive at a contradiction. Thus,  $g= E(f|\cG) \in CG^{(\cT,\cG)}_F$.

\

C4. By (C2), $CG_{G}^{(\cT, \cG)}$ is a weak$^*$ compact set and is metrizable. As argued above (\textit{e.g., Proposition~\ref{prop-CG-hemicontinuity}~(G2)}), $CG_{F_y}^{(\cT, \cG)}$ can be viewed as a compact valued correspondence from $Y$ to a compact metric space. It suffices to show that if $y_{n}\rightarrow y_0$ in $Y$ (with $y_n \ne y_0$ for each $n \ge 1$), $g_n$ weak$^*$ converges to $g_0$ in $L_p^\cG(T,X^*)$ with $f_n$ being a $\cT$-measurable and Gel$^\prime$fand integrable selection of $F_{y_n}$ and $g_n = E(f_n|\cG)$ for $n \ge 1$, then $g_0 \in CG_{F_{y_0}}^{(\cT,\cG)}$.

%For $\lambda$-almost all $t\in T$, the topological limit superior of the sequence $\{F(t, y_n)\}_{n\in \bN}$ in $(Z, d_w)$ is still the same as $w^* - \limsup_{n\to\infty} F(t, y_n)$ in $X^*$.
Since $F(t, y)$ is weak$^*$ upper hemicontinuous at $y_0$, $w^* - \limsup_{n\to\infty} F(t, y_n) \subseteq F(t, y_0)$. By Lemma~\ref{lem-distribution}, there is a subsequence of $\{f_{n}\}_{n\in\bN}$, say itself,  and a $\cT$-measurable Gel$^\prime$fand integrable selection $f$ of $F_{y_0}$ such that $\mu^{f_{n}|\cG }$ weakly converges to $\mu^{f|\cG}$. Repeating the argument in the proof of Theorem~\ref{thm-ce-bochner}~(A2) (see also the proof of (C2)), given any $\varphi$ from $(T,\cG,\lambda)$ to $X$ such that $\varphi(t)\in \{x_m\}_{m\in \bN}$ for $\lambda$-almost all $t\in T$,
$$\lim_{n\to\infty} \int_T \varphi^{**}(t,g_n(t)) \rmd\lambda  = \int_T \varphi^{**}(t,g(t)) \rmd\lambda,
$$
where $g = E(f|\cG)$. By Lemma~\ref{lem-separate_2}, $g_0 = g \in CG_{F_{y_0}}^{(\cT,\cG)}$.
\end{proof}

The following result is a straightforward corollary of Theorem~\ref{thm-ce-G}~(C4).

\begin{coro}\label{coro-CG weak upper}
Suppose that $G$ is an $\cF$-measurable, $p$-integrably bounded, and weak$^*$ compact valued correspondence from $T$ to $X^*$ for $1 < p \le \infty$. Let $\{F_n\}^\infty_{n=1}$ be a sequence of $\cF$-measurable correspondences from $T$ to $X^*$ such that $F_n(t)\subseteq G(t)$ for $\lambda$-almost all $t\in T$ and $n \ge 1$. Let $F_0 = w^*-\limsup_{n \to \infty} F_n$. Then $w^*-\limsup_{n\to\infty} CG_{F_n}^{(\cT,\cG)} \subseteq CG_{F_0}^{(\cT,\cG)}$.
\end{coro}

%%%%%%%%%%%%%%%%%%%%%%%%%%%%%%%%%%%%%%%%%%%%%%%%%%%%%%%%%%%%%%%%%%%%%%%%%%%%%%%%%%%%%%%%%%%%%%%%%%%%%%%%%%%%%%%%%%%%%%%%%%%%%%%%%%%%%%%%%%%%%%%%%%%%%%%%%%%%%%%%%%%%%%%%%%%%%%%%%%%%%%%%%%%%%%%%%%%%%%%%%%%%%%%%%%%%%%%%%%%%%%%%%%%%%%%%

\subsection{Proofs of the results in Section~\ref{sec-converse}}

To avoid repetition, we only prove the necessity results in terms of the weak and norm topologies. We shall explain how to modify the proof for the necessity results based on the weak$^*$ topology in Remark~\ref{rmk-G-necessity}.

\begin{proof}[Proof of Theorem~\ref{thm-converse}]

If $(T,\cF,\lambda)$ is purely atomic, then we are done. Suppose that $T$ can be partitioned into two disjoint parts $T_1$ and $T_2$ such that $\cF^{T_1}$ is atomless and  $\cF^{T_2}$ is atomic under $\lambda$, $T=T_1\cup T_2$, $\lambda(T_1)=1-\gamma$ with $0 \le \gamma < 1$.

Let $(I,\cI,\eta)$ be the Lebesgue unit interval. Since $(T_1,\cF^{T_1},\lambda)$ is atomless and $\cF^{T_1}$ is countably generated, by Lemma~6 in \citet{HSS2017}, there is a measure preserving mapping $\phi \colon (T_1, \cF^{T_1}, \lambda) \rightarrow ((\gamma,1], \cI_1, \eta_1)$\footnote{Let $\cI_1$ be the restriction of $\cI$ on $(\gamma,1]$, and $\eta_1$ be the Lebesgue measure on $((\gamma,1],\cI_1)$.} such that for any $E \in \cF^{T_1}$, there exists some $E' \in \cI_1$ with $\lambda( E \triangle \phi^{-1}(E'))=0$. We extend the definition of $\phi$ by letting  $\phi(t)=\frac{\gamma}{2}$ for $t\in T_2$.

For a nonnegative integer $n$ and a real number $l \in [0,1]$, the binary representation of $n$ and $l$ are given as follows:
$$n = n_0 + 2n_1 + \cdots + 2^{a-1} n_{a-1},\quad l = \frac{l_0}{2} + \frac{l_1}{2^2} + \cdots,
$$
where $n_0,n_1,\ldots,n_{a-2},l_0,l_1,\ldots \in\{0,1\}$ and $n_{a-1} = 1$. The $n$-th Walsh function $W_n$ from $(I,\cI,\eta)$ to $\{1,-1\}$ is given by
$$W_n(l)=(-1)^{n_0l_0+n_1l_1+\cdots+n_{a-1}l_{a-1}}.
$$
Let $W_0(l)=1$. Then $\{W_n(\cdot)\}_{n=0}^\infty$ is a complete orthogonal basis of $L_2(I, \bR)$, see \citet{Walsh1923}. Define $E_n =\{l\in [0,1] \colon W_n(l)=1\}$ for each integer $n \ge 0$. Since $W_0 (l)=1$ for any $l \in [0,1]$, and $\{W_n(l)\}_{n=0}^\infty$ are orthogonal, $\eta(E_n)=\frac{1}{2}$ for each $n \ge 1$.

Given an infinite-dimensional Banach space $X$, there exists a sequence of pairs $\{(x_n^*,x_n)\}_{n=0}^\infty$ such that $x_n^*\in X^*$, $x_n\in X$, $x_m^*(x_n)=0$ if $m \neq n$, and $x_n^*(x_n)=1$ for any $n$; see, for example,   \citet[Proposition 1.f.3]{LT1977}).

\

E1. Pick $\cG$ as the trivial $\sigma$-algebra; that is, $\cG = \{T,\emptyset\}$. Fix $k \geq 1$. For $j=1,\cdots, k$, define a mapping $f_j$ from $I$ to $X$ as follows:
$$
f_j(l) =
\begin{cases}
\sum_{n=0}^\infty\frac{x_{kn+j-1}}{2^n\|x_{kn+j-1}\|}W_n(\frac{l-\gamma}{1-\gamma}), & \gamma < l\leq 1 \\
0, & l\in [0,\gamma].
\end{cases}
$$
It is easy to see that $f_j$ is Bochner integrable in $X$.\footnote{For each $j$, take $s_m(l)=\sum_{n=0}^m \frac{x_{kn+j-1}}{2^n\|x_{kn+j-1}\|}W_n(\frac{l-\gamma}{1-\gamma})$. Then $\{s_m(l)\}_{m=0}^\infty$ is a sequence of Bochner integrable simple function, and $\lim_{m\to 0}\int_0^1\|f_j(l) - s_m(l)\|\rmd \eta(l)=0$.} Let $e_j = \int_I f_j(l) \rmd \eta(l)$.

%\footnote{\citet{SZ2015} had a similar, but more complex construction, which is a nonatomic game with infinite-dimensional action space and without pure strategy equilibrium. .}

Define a correspondence $F$ as follows: $F(t)=\{0, f_1\circ \phi(t), \cdots, f_k\circ \phi(t)\}$ for $t\in T$. Then $F$ is an $\cF$-measurable and weakly compact valued correspondence. We have that
\begin{align*}
\int_T f_j (\phi(t)) \rmd \lambda(t)
& = \int_{T_1} f_j (\phi(t)) \rmd \lambda(t) + \int_{T_2} f_j (\phi(t)) \rmd \lambda(t) \\
& = \int_{\gamma}^1 f_j (l) \rmd \eta(l) + \gamma \cdot f_j(\frac{\gamma}{2}) \\
& = \int_I f_j(l) \rmd \eta(l) \\
& = e_j.
\end{align*}
If $I_{F}^{\cT}$ is convex, then $\frac{1}{k+1}\sum_{j=1}^k e_j\in I_{F}^{\cT}$. Thus, there are $k+1$ disjoint $\cT$-measurable sets $S_0, S_1, \cdots, S_k \subseteq T_1$, and a $\cT$-measurable selection $f$ of $F$ such that $\int_T f \rmd \lambda= \frac{1}{k+1}\sum_{j=1}^k e_j$, $\cup_{0\leq j\leq k} S_j = T_1$ and
$$f(t)=
\begin{cases}
0,  & \text{if } t\in S_0 \cup T_2;\\
f_j(\phi(t)),  & \text{if } t\in S_j \mbox{ for } j=1,\cdots, k.
\end{cases}
$$
Due to the choice of $\phi$, there exist disjoint sets $H_0,H_1,\cdots, H_k \subseteq (\gamma, 1]$ such that $\lambda[S_j\triangle \phi^{-1}(H_j)]=0$ for $j=0,1,\cdots, k$.
Therefore,
\begin{equation*}\begin{split}\label{equa-1}
\sum_{j=1}^k \int_{H_j} f_j(l) \rmd \eta (l)
& = \sum_{j=1}^k \int_{S_j} f_j(\phi(t)) \rmd \lambda (t) = \frac{1}{k+1}\sum_{j=1}^k e_j \\
& = \frac{1}{k+1} \sum_{j=1}^k \int_{(\gamma,1]} f_j(l) \rmd\eta.
\end{split}
\end{equation*}
For each $n \geq 0$, define $\overline{E_n}=\{l\in (\gamma,1]\colon \frac{l-\gamma}{1-\gamma} \in E_n\}$ and $\overline{E_n^c}=\{l\in (\gamma,1]\colon \frac{l-\gamma}{1-\gamma} \in E_n^c\}$. For each integer $n$ and $1\leq j\leq k$, apply $x_{kn+j-1}^*$ on both sides of the equation above, we get
$$\eta\left( H_j \cap \overline{E_n} \right)-\eta\left( H_j \cap \overline{E_n^c}\right)
= \frac{1}{k+1} \left( \eta (\overline{E_n} )-\eta(\overline{E_n^c}) \right).
$$
Choosing $n=0$, $\eta(H_j)=\frac{1-\gamma}{k+1}$. In addition, for any $n$, $\eta\left( H_j \cap \overline{E_n} \right)-\eta\left( H_j \cap \overline{E_n^c}\right)$ is a constant for  all $0\leq j\leq k$.
For $n\geq 1$, since
$$0= \eta (\overline{E_n} )-\eta(\overline{E_n^c}) =\sum_{j=0}^k \eta\left( H_j \cap \overline{E_n} \right)-\eta\left( H_j \cap \overline{E_n^c}\right),
$$
we have that $\eta\left( H_j \cap \overline{E_n} \right) = \eta\left( H_j \cap \overline{E_n^c}\right)$ for $0 \le j \le k$. Then
$$\eta\left( H_j \cap \overline{E_n} \right) = \eta\left( H_j \cap \overline{E_n^c}\right) = \frac{1-\gamma}{2(k+1)}.
$$

Define $D_n=\phi^{-1}(\overline{E_n})$ for each $n$. Then
$$\lambda(D_n)= \lambda\circ\phi^{-1}(\overline{E_n})=\eta(\overline{E_n})=\frac{1-\gamma}{2}.
$$
As a result, for $n\in N$ and $0\leq j\leq k$,
\begin{align*}
\lambda(S_j\cap D_n)
& = \lambda(\phi^{-1}(H_j)\cap D_n)= \lambda\circ\phi^{-1}(H_j\cap \overline{E_n})=\eta(H_j\cap \overline{E_n}) \\
& = \frac{1-\gamma}{2(k+1)} = \frac{1}{1-\gamma}\eta(H_j)\eta(\overline{E_n})= \frac{1}{1-\gamma} \lambda(S_j)\lambda(D_n).
\end{align*}
Thus, for each $0\leq j\leq k$ and $n\in \bN$, $S_j$ is independent of $D_n$ under the probability measure $\lambda^{T_1}$. Since $\{W_n(\cdot)\}_{n=0}^\infty$ is a complete orthogonal basis of $L_2(I, \bR)$, the sequence $\{E_n\}_{n\in N}$ with null sets can generate the $\sigma$-algebra on the Lebesgue unit interval. Note that the mapping $\frac{l - \gamma}{1 - \gamma}$ is linear. It is easy to see that the sequence $\{\overline{E_n}\}_{n\in N}$ with null sets can generate the $\sigma$-algebra on $(\gamma, 1]$, which further implies that for each $0\leq j\leq k$, $S_j$ is independent of $\cF^{T_1}$ under the probability measure $\lambda^{T_1}$. By Lemma~\ref{lem-neq}, $\cT$ is nowhere equivalent to $\cF$.

\

E2/E5. Let $\cG$ be the trivial $\sigma$-algebra $\{\emptyset, T\}$ in (E2). Fix $k\geq 1$. For $1\leq j\leq k$ and $m\geq 1$, let $f_j^m$ be a mapping from $I$ to the space spanned by $\{x_0,\cdots,x_{km+k-1}\}$ as follows:
$$
f_j^m(l)=\begin{cases}
\sum_{n=0}^m\frac{x_{kn+j-1}}{2^n\|x_{kn+j-1}\|}W_n(\frac{l-\gamma}{1-\gamma}), & \gamma < l\leq 1; \\
0, & l\in [0,\gamma].
\end{cases}
$$
Then $e_j = \int_I f_j^m(l) \rmd\eta$ for $1\leq j \leq k$. Let $F^m=\{0,f^m_1\circ \phi,\cdots, f^m_k\circ \phi\}$.  By Lyapunov theorem (see Proposition~5 in \citet[p.62]{Hd1974}), $I_{F^m}^\cT$ is convex. Thus, we have $k+1$ disjoint   sets $S^m_0, S^m_1, \cdots, S^m_k \subseteq T_1$ and a $\cT$-measurable selection $f^m$ of $F^m$  such that $\int_T f^m \rmd\lambda= \frac{1}{k+1}\sum_{j=1}^k e_j$, where
$$f^m(t)=
\begin{cases}
0,  & \text{if } t\in S^m_0 \cup T_2;\\
f_j^m(\phi(t)),  & \text{if } t\in S^m_j \mbox{ for } j=1,\cdots, k.
\end{cases}$$

Consider the correspondence $F$ constructed in the proof of (E1). Since $F$ is norm (and also weakly) compact valued, due to the condition, $I_F^\cT$ is norm (and also weakly) compact.
Define a mapping $g^m$ as follow:
$$g^m(t)=
\begin{cases}
0,  & \text{if } t\in S^m_0 \cup T_2;\\
f_j(\phi(t)),  & \text{if } t\in S^m_j \mbox{ for } j=1,\cdots, k.
\end{cases}$$
Then $\int_{T} g^m \rmd\lambda \in I_F^\cT$. It is clear that $\int_{T} g^m \rmd\lambda$ both norm and weakly converges to $\frac{1}{k+1}\sum_{j=1}^k e_j$, thus $\frac{1}{k+1}\sum_{j=1}^k e_j\in I_F^\cT$. By the proof of (E1), we are done.

\

E3. Consider the correspondence $F$ in the proof of (E1). Since $I_{F}^{\cT}=I_{\overline{co}F}^{\cT}$, we have that
$$\frac{1}{k+1}\sum_{j=1}^k e_j = \frac{1}{k+1}\sum_{j=1}^k \int_I f_j(l) \rmd\eta = \frac{1}{k+1}\sum_{j=1}^k \int_T f_j(\phi(t)) \rmd\lambda \in I_{\overline{co}F}^{\cT} = I_{F}^{\cT}.
$$
Repeating the argument in the proof of (E1), we are done.

\

E4/E6. Let $\cG$ be the trivial $\sigma$-algebra $\{\emptyset, T\}$ in (E4). Fix $k\geq 1$ and $m\geq 1$. Recall the correspondence $F^m$ defined in the proof of (E2/E5). Let $Y$ be the space $\{0,1,\frac{1}{2},\cdots, \frac{1}{m},\cdots\}$ endowed with the usual metric. Let $\Phi$ be the correspondence from $T \times Y$ to $X$ such that for each $t\in T$, $\Phi(t, \frac{1}{m})= F_m(t)$ for $m\geq 1$ and $\Phi(t,0) = F(t)$, where $F(t)$ has been defined in the proof of (E1). It is easy to see that the correspondence $\Phi$ satisfies the conditions specified in (E4/E6).

We have shown in the proof of (E2/E5) that for each $m$,
$$\frac{1}{k+1}\sum_{j=1}^k e_j \in I_{F_m}^\cT = I_{\Phi_{\frac{1}{m}}}^\cT.
$$
Thus, $\frac{1}{k+1}\sum_{j=1}^k e_j \in I_{\Phi_0}^\cT = I_{F}^\cT$. The proof can then be completed by repeating the argument in the proof of (E1).
\end{proof}

\begin{rmk}\label{rmk-G-necessity}
To prove the necessity results based on the weak$^*$ topology, we shall make the following changes.
\begin{itemize}
\item In the proof of (E1), define
$$f_j(l) =
\begin{cases}
\sum_{n=0}^\infty\frac{x^*_{kn+j-1}}{2^n\|x^*_{kn+j-1}\|}W_n(\frac{l-\gamma}{1-\gamma}), & \gamma < l\leq 1 \\
0, & l\in [0,\gamma].
\end{cases}
$$

\item In the proof of (E2/E5), define
$$f_j^m(l) =
\begin{cases}
\sum_{n=0}^m\frac{x_{kn+j-1}}{2^n\|x_{kn+j-1}\|}W_n(\frac{l-\gamma}{1-\gamma}), & \gamma < l\leq 1; \\
0, & l\in [0,\gamma].
\end{cases}
$$

\item In the proofs of (E3) and (E4/E6), no change.
\end{itemize}
\end{rmk}

%%%%%%%%%%%%%%%%%%%%%%%%%%%%%%%%%%%%%%%%%%%%%%%%%%%%%%%%%%%%%%%%%%%%%%%%%%%%%%%%%%%%%%%%%%%%%%%%%%%%%%%%%%%%%%%%%%%%%%%%%%%%%%%%%%%%%%%%%%%%%%%%%%%%%%%%%%%%%%%%%%%%%%%%%%%%%%%%%%%%%%%%%%%%%%%%%%%%%%%%%%%%%%%%%%%%%%%%%%%%%%%%%%%%%%%%

\subsection{Proofs of the results in Section~\ref{sec-game}}\label{subsec-game proof}

We first show the easy direction of Theorem \ref{thm-existence} that (1) $\Rightarrow$ (2).

\paragraph{(1) $\Rightarrow$ (2).}
We first focus on the case with Bochner integrable strategy profile. Given any $\gamma \in \Gamma$, consider the best-response correspondence $F$ from $T \times \Gamma$ to $A$,
$$F(t, \gamma) = \mbox{argmax}_{a\in A} G_t(a, \gamma),
$$
and the correspondence $H$ from $\Gamma$ to $\Gamma$,
$$H(\gamma) = \big\{ E(f(t, \gamma) | \cF) \colon f(\cdot, \gamma) \text{ is a } \cT \text{-measurable selection of } F(\cdot, \gamma) \big\}.
$$
We shall show that $H$ is nonempty and convex valued, and weakly upper hemicontinuous.

For any $\gamma$, $G_{t}(a, \gamma)$ is $\cF$-measurable in $t$ and weakly continuous in $a$. By Lemma~III.39 in \citet{CV1977}, $G_{t}(a, \gamma)$ is $\cF \times \cB(A)$-measurable on $T \times A$. By Berge's maximum theorem (see Theorem~17.31 in \citet{AB2006}), $F(t, \cdot)$ is weakly compact valued and upper hemicontinuous in $\gamma$. In addition, $F$ is $\cF \otimes \cB(\Gamma)$-measurable and admits a measurable selection due to Theorem~18.19 in \citet{AB2006}. Thus, there exists an $\cF$-measurable selection of the correspondence $F_{\gamma}$ on $T$, where $F_{\gamma}(\cdot)=F(\cdot, \gamma)$. It further implies that $H(\gamma)$ is nonempty.

Let $F'$ be a constant correspondence from $T$ to $A$, where $F'(t) = A$. Then $F'$ is an integrably bounded and weakly compact valued correspondence. Besides, $F(t, \gamma) \subseteq A = F'(t)$ for each $(t, \gamma)$. Applying Theorem~\ref{thm-ce-bochner}, $H(\gamma)$ is convex-valued and weakly upper hemicontinuous. By Proposition~\ref{prop-CI-hemicontinuity}~(F1), $\Gamma$ is convex and weakly compact. By the Fan-Glicksberg fixed point theorem (see Corollary~17.55 in \citet{AB2006}), the correspondence $H$ admits a fixed point $\gamma \in \Gamma$. By the definition of $H$, there exists a $\cT$-measurable mapping $g \colon T \to A$ such that $g(t) \in F(t, \gamma)$ and $\gamma = E(g | \cF)$. It is clear that $g$ is a Nash equilibrium of the game.

By applying Theorem~\ref{thm-ce-G}, the argument above still works in the case with Gel$^\prime$fand integrable strategy profile.

\

Next, we are going to prove the opposite direction (2) $\Rightarrow$ (1).

\paragraph{(2) $\Rightarrow$ (1).}
If $(T,\cF,\lambda)$ is purely atomic, then we are done. Suppose that $T$ can be partitioned into two disjoint parts $T_1$ and $T_2$ such that $\cF^{T_1}$ is atomless and  $\cF^{T_2}$ is atomic under $\lambda$, $T=T_1\cup T_2$, $\lambda(T_1) = 1 - \gamma$ with $0 \le \gamma < 1$. Let $(I,\cI,\eta)$ be the Lebesgue unit interval. Since $(T_1,\cF^{T_1},\lambda)$ is atomless and $\cF^{T_1}$ is countably generated, by Lemma~6 in \citet{HSS2017}, there is a measure preserving mapping $\phi \colon (T_1, \cF^{T_1}, \lambda) \rightarrow ((\gamma,1], \cI_1, \eta_1)$ such that for any $E \in \cF^{T_1}$, there exists some $E' \in \cI_1$ with $\lambda( E \triangle \phi^{-1}(E'))=0$. For $t\in T_2$, let $\phi(t)=\frac{\gamma}{2}$.

By Lemma~\ref{lem-neq}, it is sufficient to show that $\cF^{T_1}$ admits an asymptotic independent supplement in $\cT^{T_1}$ under the probability measure $\lambda^{T_1}$. In particular, we only need to prove the following claim.

\begin{claim}\label{claim-asymptotic}
Suppose that Condition~(1) in Theorem \ref{thm-existence} hold. For any integer $k \ge 2$, there is a $\cT$-measurable partition $\{P_i\}_{i=0}^{k}$ of $T_1$ such that $\lambda(P_i) = \frac{1 - \gamma}{k+1}$ for $i = 0,1,\cdots, k$, $P_i\cap P_j = \emptyset$ for $i \neq j$, and $\{P_i\}_{i=0}^{k}$ are independent of $\cF^{T_1}$ under the probability measure $\lambda^{T_1}$.
\end{claim}

\begin{proof}
We first consider the case with Bochner integrable strategy profile.

Denote $\alpha$ as the $(k+1)$-th root of $1$, i.e., $\alpha^j \neq 1$ for $j = 1, \ldots, k$, and $\alpha^{k+1}=1$. Recall the construction of the Walsh system in the proof of Theorem~\ref{thm-converse}. Given a nonnegative integer $n$ and a real number $l \in I$, the $n$-th Walsh function $W_n$ from $I$ to $\{1,-1\}$ is given by
$$W_n(l)=(-1)^{n_0l_0+n_1l_1+\cdots+n_{a-1}l_{a-1}},
$$
where
$$n = n_0 + 2n_1 + \cdots + 2^{a-1} n_{a-1} \quad \mbox{ and } \quad l = \frac{l_0}{2} + \frac{l_1}{2^2} + \cdots,
$$
with $n_0,n_1,\ldots,n_{a-2},l_0,l_1,\ldots \in\{0,1\}$ and $n_{a-1} = 1$. Let $W_0(l) = 1$. Then $\{W_n(\cdot)\}_{n=0}^\infty$ is a complete orthogonal basis of $L_2(I, \bR)$. Define $E_n =\{l\in I \colon W_n(l)=1\}$ for each integer $n \ge 0$. We have that $\eta(E_n) = \frac{1}{2}$ for $n \ge 1$, and $\eta(E_0) = 1$.

For the infinite-dimensional Banach space $X$, there exists a sequence of pairs $\{(x_n^*,x_n)\}_{n=0}^\infty$ such that $x_n^*\in X^*$, $x_n\in X$, $x_m^*(x_n)=0$ if $m\neq n$, and $x_n^*(x_n)=1$ for all $n$ (see Proposition 1.f.3 in \citet{LT1977}).

Define $\psi_1$, $\psi_2$, $\cdots$, $\psi_k\colon I\to X$ as follows:
$$\psi_1(l) =
\begin{cases}
\sum_{n=0}^\infty\frac{x_{kn}}{2^n\|x_{kn}\|}W_n(\frac{l-\gamma}{1-\gamma}), & \gamma < l\leq 1 \\
0, & l\in [0,\gamma].
\end{cases}
$$
$$\psi_2(l) =
\begin{cases}
\sum_{n=0}^\infty\frac{x_{kn+1}}{2^n\|x_{kn+1}\|}W_n(\frac{l-\gamma}{1-\gamma}), & \gamma < l\leq 1 \\
0, & l\in [0,\gamma].
\end{cases}
$$
$$\cdots$$
$$\psi_k(l) =
\begin{cases}
\sum_{n=0}^\infty\frac{x_{kn+k-1}}{2^n\|x_{kn+k-1}\|}W_n(\frac{l-\gamma}{1-\gamma}), & \gamma < l\leq 1 \\
0, & l\in [0,\gamma].
\end{cases}
$$
Denote $e_i=\int_I \psi_i(l) \rmd \eta(l)$ for $i = 1, 2,\cdots, k$. Let $\Psi(l)=\{0,\psi_1(l),\cdots,\psi_k(l)\}$, where $0$ is the zero function from $I$ to $X$.

%that is,
%$$\overline{co}(A) \supseteq \left\{ \int_I g\rmd\eta \colon \text{$g$ is a Bochner integrable function from $I$ to $A$} \right\}.
%$$

Now we are ready to construct a game. The player space is $(T, \cT, \lambda)$. The action space $A$ is a convex and weakly compact set in $X$ that contains $\{\Psi(l)\}_{l \in I}$.\footnote{Such a set $A$ exists since $\{\Psi(l)\}_{l \in I}$ is bounded.} Hence, $\overline{co}(A) = A$. By Mazur's theorem, the set $A$ is both weakly compact and norm compact. Let $M = \sup\{\|x\|\colon x\in A\}$. Note that $\left\|e_1\right\|=\left\|(1 - \gamma) \frac{x_0}{\|x_0\|}\right\| = 1 - \gamma$, which implies that $M \ge 1 - \gamma$. For each player~$t$, the payoff $G(t)$ only depends on the player's action and the societal aggregate actions. Since $\|b\|\le M$ for any $b \in A$, and $\|e_i\|\le M$ for $1\leq i\leq k$,
$$\|b-\frac{e_1+\cdots+e_k}{k+1}\|\le 2M.
$$
Denote the distance in $X$ by $d(\cdot,\cdot)$. That is, $d(x, y) = \|x - y\|$ for any $x, y \in X$. We have that for any $b \in A$
$$\beta d\left(b,\frac{e_1+\cdots+e_k}{k+1}\right) \equiv \frac{1 - \gamma}{4M} \left\|b - \frac{e_1+\cdots+e_k}{k+1} \right\| \le \frac{1 - \gamma}{2} < \frac{1}{2},
$$
where $\beta = \frac{1 - \gamma}{4M} \le \frac{1}{4}$.

Let $h$ be a mapping from $[0,1] \times A \times A^k \times R_+$ to $R_+$ as follows:
$$h(l, a, x_1, \cdots, x_k, \theta) =
$$
$$\begin{cases}
0  & \text{if }\theta = 0 \mbox{ or } l \in [0, \gamma], \\
\theta \left|\sin\frac{l - \gamma}{\theta} \pi \right| \cdot (\|a\|+\left|1-\alpha^{\left[\frac{l - \gamma}{\theta}\right]}\right|)
\\ \cdot \prod_{i=1}^{k}\left(\|a-\frac{1}{k+1}x_i-\frac{1}{k+1}\sum_{j=1}^{k}x_j\|
+\left|\alpha^{i}-\alpha^{\left[\frac{l - \gamma}{\theta}\right]}\right|\right)  & \text{otherwise },
\end{cases}
$$
where $[\frac{l - \gamma}{\theta}]$ is the integer part of the real number $\frac{l - \gamma}{\theta}$. The payoff function $G \colon T\to \cU$ is defined as follows: for any player $t \in T$, given the action $a \in A$ and the societal aggregate $b \in A$,
\begin{align*}
G(t)\left(a,b\right)
& = -h\left(\phi(t),a,\psi_1(\phi(t)),\cdots,\psi_k(\phi(t)),\beta d\left(b,\frac{e_1+\cdots+e_k}{k+1}\right) \right) \\
& \quad -\|a\| \cdot \prod_{i=1}^{k} \bigg\|a-\frac{1}{k+1}\psi_i(\phi(t)) - \frac{1}{k+1} \sum_{j=1}^{k}\psi_j(\phi(t)) \bigg\|.
\end{align*}
The game $G$ is $\cF$-measurable as $\phi$ is $\cF$-measurable, and norm continuous in $(a, b)$ (and hence weakly continuous in $(a, b)$). If player $t$ takes action $a$ and all other players follow the strategy $f$, the payoff for player $t$ is $G(t)(a,\int_T f \rmd \lambda)$.

By Condition~(1), there exists a $\cT$-measurable strategy $g$ that is a Nash equilibrium of the game $G$. For player $t \in T_2$, it is clear that the unique best response is $a = 0$, which implies that $g(t) = 0$. Thus, we only need to identify the best response for player $t \in T_1$. There are two possible cases to study.

\

Case 1: Suppose that $\int_{T} g(t) \rmd \lambda \neq \frac{e_1 + \cdots + e_k}{k+1}$.

Let $$d_0=\beta d\left(\int_{T}g(t)\rmd \lambda,\frac{e_1+\cdots+e_k}{k+1}\right).
$$
Then $0 < d_0 \leq \frac{1}{2}$. Partition the interval $I_1 = (\gamma, 1]$ by the length $d_0$, we get $\{I_1 \cap \left[\gamma + nd_0, \gamma + (n+1)d_0\right] \}_{n \ge 0}$. Fix some player $t \in T_1$ with $\phi(t) \in \left(\gamma + (kj + i)d_0, \gamma + (kj + i + 1)d_0\right)$ for $j \ge 0$ and $i = 0, 1, \ldots, k$. We consider the best response of player~$t$.

If $i=0$, then
$$G_t\left(a,\int_{T}g\rmd \lambda\right)=
$$
$$ - d_0 \left|\sin\frac{\phi(t) - \gamma}{d_0}\pi\right| \cdot\|a\| \cdot \prod_{m=1}^{k}\left(\left\|a-\frac{1}{k+1} \psi_m(\phi(t)) -\frac{1}{k+1}\sum_{j=1}^{k}\psi_j(\phi(t)) \right\| + \left|\alpha^{m} - 1 \right|\right)
$$
$$ - \|a\| \cdot \prod_{m=1}^{k} \bigg\|a - \frac{1}{k+1}\psi_m(\phi(t)) - \frac{1}{k+1}\sum_{j=1}^{k}\psi_j(\phi(t)) \bigg\|.
$$
Since $G_t$ is always negative and its value is $0$ if and only if $a = 0$, the best response for player $t$ is $a=0$. That is, $g(t)=0$ in this case.

Following the same argument, one could show that when $i = 1, \ldots, k$, the best response is
$$g(t) = \frac{1}{k + 1}\psi_i(\phi(t)) + \frac{1}{k+1} \sum_{j=1}^{k} \psi_j(\phi(t)).
$$
For $i = 1, \ldots, k$, let $q_i(l) = i$ if $l \in \left(\gamma + (kj + i)d_0, \gamma + (kj + i + 1)d_0\right)$ for some $j \ge 0$, and $0$ otherwise. Then for $t \in T_1$,
$$g(t) = \frac{1}{k+1} \sum_{j=1}^{k} \left[ q_j(\phi(t)) \cdot \psi_j(\phi(t)) \right] + \frac{1}{k+1} \sum_{j=1}^{k}q_j(\phi(t)) \cdot \sum_{j=1}^{k}\psi_j(\phi(t)).
$$

%In sum, for almost all $t\in T$,
%$$g(t)
%=\begin{cases}
%0,  & \text{if } \phi(t)\in\left(3kd_0,(3k+1)d_0\right) \text{for some } k\in N; \\
%\frac{2}{3}\psi_1(\phi(t))+\frac{1}{3}\psi_2(\phi(t)),  & \text{if } \phi(t)\in\left((3k+1)d_0,(3k+2)d_0\right) \text{for some } k\in N;\\
%\frac{1}{3}\psi_1(\phi(t))+\frac{2}{3}\psi_2(\phi(t)),  & \text{if } \phi(t)\in\left((3k+2)d_0,3(k+1))d_0\right) \text{for some } k\in N.
%\end{cases}$$

We have
\begin{align}
& \quad d(\int_T g(t) \rmd \lambda(t), \frac{e_1 + \ldots + e_k}{k+1}) \nonumber \\
& = \bigg\| \frac{1}{k+1} \sum_{j=1}^{k} \int_{T_1} q_j(\phi(t)) \cdot \psi_j(\phi(t)) \rmd \lambda(t) + \frac{1}{k+1} \int_{T_1} \sum_{j=1}^{k}q_j(\phi(t)) \cdot \sum_{j=1}^{k}\psi_j(\phi(t)) \rmd \lambda(t) \nonumber \\
& - \frac{1}{k + 1} \int_{\gamma}^{1} \psi_1 \rmd \eta - \cdots - \frac{1}{k + 1} \int_{\gamma}^{1} \psi_k \rmd \eta \bigg\| \nonumber \\
& \le \frac{1}{k + 1} \bigg\| \int_{\gamma}^{1} q_1 \cdot \psi_1 \rmd \eta + \int_{\gamma}^{1} \psi_1 \cdot \sum_{j=1}^{k} q_j \rmd \eta - \int_{\gamma}^{1} \psi_1 \rmd \eta \bigg\| \nonumber \\
& + \cdots \nonumber \\
& + \frac{1}{k + 1} \bigg\| \int_{\gamma}^{1} q_k \cdot \psi_k \rmd \eta + \int_{\gamma}^{1} \psi_k \cdot \sum_{j=1}^{k}q_j \rmd \eta - \int_{\gamma}^{1} \psi_k \rmd \eta \bigg\|.
\end{align}

To estimate (1), we need the following preparatory lemma.\footnote{This lemma is analogous to Lemmas~2, 3, 4 in \citet{SZ2015}, we provide the detailed proof here for completeness.}

\begin{lem}\label{lem-inequ}
For any $i =1, \ldots, k$,
$$\sum_{n=0}^\infty \frac{1}{2^n} \left| \int_{\gamma}^{1} \left[ q_i (l) + \sum_{j=1}^{k} q_j(l) - 1 \right] W_n(\frac{l-\gamma}{1-\gamma}) \rmd\eta(l) \right| < 4 d_0.
$$
\end{lem}

\begin{proof}
For any integer $n \ge 0$ and $i =1, \ldots, k$, we first estimate the absolute bound for $\int_{\gamma}^{1} \left[ q_i (l) + \sum_{j=1}^{k} q_j(l) - 1 \right] W_n(\frac{l-\gamma}{1-\gamma}) \rmd\eta(l)$. It is obvious that
$$\left| \int_{\gamma}^{1} \left[ q_i (l) + \sum_{j=1}^{k} q_j(l) - 1 \right] W_n(\frac{l-\gamma}{1-\gamma}) \rmd\eta(l) \right| \le 1.
$$
Let the binary representation of $n$ be
$$n = n_0 + 2n_1 + 2^2n_2 + \cdots + 2^{m-1}n_{m-1},
$$
where $n_{m-1} =1$ and $n_0,n_1,\ldots, n_{m-2} \in \{0, 1\}$. Given an integer $s \in \{1,2,\ldots,2^m\}$,
for any $\frac{l-\gamma}{1-\gamma} \in \left( \frac{s-1}{2^m}, \frac{s}{2^m} \right)$, the first $m$ components in the binary representation of $\frac{l-\gamma}{1-\gamma}$ are fixed. By its definition, $W_n(\frac{l-\gamma}{1-\gamma}) = (-1)^{c_s}$ for some constant integer $c_s$ when $\frac{l-\gamma}{1-\gamma} \in \left( \frac{s-1}{2^m}, \frac{s}{2^m} \right)$. We have that
\begin{align*}
& \quad \left| \int_{\gamma}^{1} \left[ q_i (l) + \sum_{j=1}^{k} q_j(l) - 1 \right] W_n(\frac{l-\gamma}{1-\gamma}) \rmd \eta(l) \right| \\
& \le \sum_{s=1}^{2^m} \left| \int_{\gamma + (1 - \gamma) \frac{s-1}{2^m}}^{\gamma + (1 - \gamma) \frac{s}{2^m}} \left[ q_i (l) + \sum_{j=1}^{k} q_j(l) - 1 \right] W_n(\frac{l-\gamma}{1-\gamma}) \rmd \eta(l) \right| \\
& = \sum_{s=1}^{2^m} \left| \int_{\gamma + (1 - \gamma) \frac{s-1}{2^m}}^{\gamma + (1 - \gamma) \frac{s}{2^m}} \left[ q_i (l) + \sum_{j=1}^{k} q_j(l) - 1 \right] (-1)^{c_s} \rmd \eta(l) \right| \\
& \le 2^m d_0 \\
& \le 2n d_0.
\end{align*}
The second inequality holds since for any $s = 1, \ldots, 2^m$,
$$\left| \int_{\gamma + (1 - \gamma) \frac{s-1}{2^m}}^{\gamma + (1 - \gamma) \frac{s}{2^m}} \left[ q_i (l) + \sum_{j=1}^{k} q_j(l) - 1 \right] (-1)^{c_s} \rmd \eta(l) \right| \le d_0.
$$
The last inequality is due to the binary representation of $n$. Thus, we have
$$ \left| \int_{\gamma}^{1} \left[ q_i (l) + \sum_{j=1}^{k} q_j(l) - 1 \right] W_n(\frac{l-\gamma}{1-\gamma}) \rmd\eta(l) \right| \le \min\{1, 2n d_0 \}.
$$

Let $\kappa$ be the integer such that $2\kappa d_0 \le 1$ and $2(\kappa+1) d_0 > 1$. Then
\begin{align*}
& \quad \sum_{n=0}^\infty \frac{1}{2^n} \left| \int_{\gamma}^{1} \left[ q_i (l) + \sum_{j=1}^{k} q_j(l) - 1 \right] W_n(\frac{l-\gamma}{1-\gamma}) \rmd\eta(l) \right| \\
& < \sum_{n=0}^\infty \frac{1}{2^n} \min\{1, 2n d_0 \} \\
& = \sum_{n = 0}^{\kappa} \frac{2n d_0}{2^n} + \sum_{n=\kappa+1}^\infty \frac{1}{2^n} \\
& = 2(2 - \frac{\kappa+2}{2^{\kappa}}) d_0 + \frac{1}{2^{\kappa}} \\
& < 4d_0.
\end{align*}
This completes the proof.
\end{proof}

Now we are ready to go back to estimate (1). For $i = 1, \ldots, k$,
\begin{align*}
& \quad \frac{1}{k + 1} \bigg\| \int_{\gamma}^{1} q_i \cdot \psi_i \rmd \eta + \int_{\gamma}^{1} \psi_i \cdot \sum_{j=1}^{k} q_j \rmd \eta - \int_{\gamma}^{1} \psi_i \rmd \eta \bigg\| \\
&= \frac{1}{k + 1} \left\| \int_{\gamma}^{1} \left[ q_i(l) + \sum_{j=1}^{k} q_j(l) - 1 \right] \sum_{n=0}^\infty \frac{x_{kn+i-1}}{2^n\|x_{kn+i-1}\|} W_n(\frac{l-\gamma}{1-\gamma}) \rmd \eta(l) \right\| \\
& \le \frac{1}{k + 1} \sum_{n=0}^\infty \left\| \int_{\gamma}^{1} \left[ q_i(l) + \sum_{j=1}^{k} q_j(l) - 1 \right] \frac{x_{kn+i-1}}{2^n\|x_{kn+i-1}\|} W_n(\frac{l-\gamma}{1-\gamma}) \rmd \eta(l) \right\| \\
& \le \frac{1}{k + 1} \sum_{n=0}^\infty \frac{1}{2^n} \left| \int_{\gamma}^{1} \left[ q_i(l) + \sum_{j=1}^{k} q_j(l) - 1 \right] W_n(\frac{l-\gamma}{1-\gamma}) \rmd \eta(l) \right| \\
& \le \frac{4}{k + 1} d_0.
\end{align*}
The last inequality holds because of Lemma~\ref{lem-inequ}. Then we have $(1) \le \frac{4k}{k + 1} d_0$, which implies that $d_0 \le \frac{4k}{k + 1} \beta d_0$. Thus, $\beta \ge \frac{k + 1}{4k} > \frac{1}{4}$. This is a contradiction since $\beta = \frac{1 - \gamma}{4M} \le \frac{1}{4}$.

\

Case 2: Suppose that $\int_{T}g(t)\rmd \lambda = \frac{e_1+\cdots+e_k}{k+1}$, therefore,
$$G(t)\left(a,\int_{L}g(t)\rmd \lambda\right) = -\|a\| \cdot \prod_{i=1}^{k} \bigg\|a - \frac{1}{k+1}\psi_i(\phi(t)) - \frac{1}{k+1} \sum_{j=1}^{k}\psi_j(\phi(t)) \bigg\|.
$$
The best response correspondence is hence
$$\Psi_1(t) = \{0\} \cup \left\{ \frac{1}{k+1}\psi_i(\phi(t)) + \frac{1}{k+1} \sum_{j=1}^{k}\psi_j(\phi(t)) \right\}_{1 \le i \le k}
$$
for $t \in T_1$. The equilibrium strategy $g$ must be a $\cT$-measurable selection of the correspondence $\Psi_1$ on $T_1$. There exists a partition $\{P_0, P_1, \ldots, P_k\}$ of $T_1$ such that
$$g(t) =
\begin{cases}
0,  & \text{if } t\in P_0; \\
\frac{1}{k+1}\psi_1(\phi(t)) + \frac{1}{k+1} \sum_{j=1}^{k}\psi_j(\phi(t)),  & \text{if } t\in P_1; \\
\ldots, & \ldots; \\
\frac{1}{k+1}\psi_k(\phi(t)) + \frac{1}{k+1} \sum_{j=1}^{k}\psi_j(\phi(t)),  & \text{if } t\in P_k.\\
\end{cases}
$$
By the definition of $\phi$, there exist $k+1$ subsets $\tilde{P}_0, \tilde{P}_1, \ldots, \tilde{P}_k \subseteq (\gamma, 1]$ such that $\lambda[P_j \triangle \phi^{-1}(\tilde{P}_j)] = 0$ for $j = 0, 1, \ldots, k$.

Since $\int_{T}g(t)\rmd \lambda = \frac{e_1+\cdots+e_k}{k+1}$,
\begin{align*}
& \quad \frac{e_1 + \ldots + e_k}{k+1} \\
& = \sum_{1 \le i \le k} \int_{P_i} \left[ \frac{1}{k+1}\psi_i(\phi(t)) + \frac{1}{k+1} \sum_{j=1}^{k}\psi_j(\phi(t)) \right] \rmd \lambda(t) \\
& = \sum_{1 \le i \le k} \int_{\tilde{P}_i} \left[ \frac{1}{k+1}\psi_i(l) + \frac{1}{k+1} \sum_{j=1}^{k}\psi_j(l) \right] \rmd \eta(l).
\end{align*}
By the definition of $\{e_i\}_{1 \le i \le k}$,
$$\sum_{1 \le i \le k} \int_{\eta}^{1}  \psi_i(l) \rmd \eta(l) = \sum_{1 \le i \le k} \int_{\tilde{P}_i} \left[ \psi_i(l) + \sum_{j=1}^{k}\psi_j(l) \right] \rmd \eta(l).
$$

For each $n \ge 0$, let $\tilde{E}_n =\{l \in (\gamma, 1] \colon \frac{l - \gamma}{1 - \gamma} \in E_n \}$ and $\overline{E}_n =\{l \in (\gamma, 1] \colon \frac{l - \gamma}{1 - \gamma} \notin E_n \}$. Then $W_n(\frac{l - \gamma}{1 - \gamma}) = 1$ for $l \in \tilde{E}_n$, and $-1$ for $l \in \overline{E}_n$. In addition, since the transformation $\frac{l - \gamma}{1 - \gamma}$ is linear and $\eta(E_n) = \eta(E^c_n)$ for $n \ge 1$, $\eta(\tilde{E}_n) = \eta(\overline{E}_n)$.

For $1 \le i \le k$ and $n \ge 0$, applying $x_{kn+i}^*$ on both sides of the above equality, we get
$$\eta \left( \tilde{P}_0 \cap \tilde{E}_n \right) - \eta\left(\tilde{P}_0 \cap \overline{E}_n \right) = \eta\left(\tilde{P}_i \cap \tilde{E}_n \right)-\eta\left(\tilde{P}_i \cap \overline{E}_n \right).
$$
Pick $n=0$. As $\tilde{E}_0 = (\gamma, 1]$ and $\overline{E}_0 = \emptyset$, we get
$$\eta(\tilde{P}_0) = \eta(\tilde{P}_i) = \frac{1 - \gamma}{k + 1}.
$$
As a result, for $i = 0, \ldots, k$,
$$\lambda(P_i) = \lambda(\phi^{-1}(\tilde{P}_i)) = \eta(\tilde{P}_i) = \frac{1 - \gamma}{k + 1}.
$$
For $n \ge 1$,
$$0 = \eta(\tilde{E}_n) - \eta(\overline{E}_n) = \sum_{i=0}^k
\eta\left(\tilde{P}_i \cap \tilde{E}_n \right)-\eta\left(\tilde{P}_i \cap \overline{E}_n \right),
$$
which implies that for $i =0, \ldots, k$,
$$\eta\left(\tilde{P}_i \cap \tilde{E}_n \right) - \eta\left(\tilde{P}_i \cap \overline{E}_n \right) = 0.
$$
That is,
$$\eta\left(\tilde{P}_i \cap \tilde{E}_n \right) = \frac{1}{2} \eta\left( \tilde{P}_i \right) = \frac{1 - \gamma}{2(k + 1)}.
$$
Let $D_n = \phi^{-1}(\tilde{E}_n)$ for each $n$. Then $\lambda(D_n) = \lambda\circ\phi^{-1}(\tilde{E}_n) = \eta(\tilde{E}_n) = \frac{1 - \gamma}{2}$. Thus,
\begin{align*}
\lambda(P_i \cap D_n)
& = \lambda(\phi^{-1}(\tilde{P}_i) \cap D_n) = \lambda \circ\phi^{-1}(\tilde{P}_i \cap \tilde{E}_n) = \eta(\tilde{P}_i \cap \tilde{E}_n) \\
& = \frac{1 - \gamma}{2(k + 1)} = \frac{1}{1 - \gamma} \cdot \frac{1 - \gamma}{k + 1} \cdot  \frac{1 - \gamma}{2} \\
& = \frac{1}{1 - \gamma} \eta(\tilde{P}_i) \eta(\tilde{E}_n) = \frac{1}{1 - \gamma} \lambda(P_i)\lambda(D_n).
\end{align*}
Thus, for each $n \ge 1$ and $i =0, \ldots, k$, $P_i$ is independent of $D_n$ under the probability measure $\lambda^{T_1}$. Repeat the same argument in the last paragraph of the proof of Theorem~\ref{thm-converse}~(E1), $P_i$ is independent of $\cF^{T_1}$ under the probability measure $\lambda^{T_1}$. This completes the proof for the case with Bochner integrable strategy profile.

\

For the case with Gel$^\prime$fand integrable strategy profile, recall that $X$ is the dual space of $Y$, where $Y$ is an infinite-dimensional separable Banach space. There exists a sequence of pairs $\{(x_n, y_n)\}_{n=0}^\infty$ such that $x_n \in X$, $y_n \in Y$, $x_m(y_n) = 0$ if $m\neq n$, and $x_n(y_n) = 1$ for all $n$. The only difference in the construction of the game is that the action space $A$ is a convex and weak$^*$ compact set in $X$. Then $w^*-\overline{co}(A) = A$. The payoff function $G_t(a, b)$ is weak$^*$ continuous in $(a, b)$ as it is norm continuous. The discussion on Case~1 remains the same as above. In the discussion on Case~2, one applies $y_{kn+i}$ rather than $x_{kn+i}^*$. Then the argument still works.
\end{proof}

\begin{rmk}\label{rmk-game-norm proof}
Now we describe how to modify the proof of Theorem~\ref{thm-existence} to cover the case based on norm topology (Remark~\ref{rmk-game-norm}).

For the direction that (1) $\Rightarrow$ (2), one needs to revise $F$ and $H$ as follows: given any $x \in \overline{\mbox{co}}(A)$, consider the best-response correspondence $F$ from $T \times \overline{\mbox{co}}(A)$ to $A$,
$$F(t, x) = \mbox{argmax}_{a\in A} G_t(a, x),
$$
and the correspondence $H$ from $\overline{\mbox{co}}(A)$ to $\overline{\mbox{co}}(A)$,
$$H(x) = \bigg\{ \int_{T} f(t, x) \rmd \lambda(t) \colon f(\cdot, x) \text{ is a } \cT \text{-measurable selection of } F(\cdot, x) \bigg\}.
$$
Then one can follow the same argument as that in the proof of Theorem~\ref{thm-existence} and adopt Proposition~\ref{prop-norm integral} to identify a fixed point $x \in \overline{\mbox{co}}(A)$ for $H$, and a $\cT$-measurable mapping $g \colon T \to A$ such that $g(t) \in F(t, x)$ and $x = \int g \rmd \lambda$. It is clear that $g$ is a Nash equilibrium of the game.

For the direction that (2) $\Rightarrow$ (1), just note that the constructed action set $A$ is also norm compact, and the payoff function $G_t(\cdot, \cdot)$ is indeed norm continuous. Hence, the proof of Theorem~\ref{thm-existence} is also applicable.
\end{rmk}

%\newpage

{\small
%\singlespacing

\end{document}